\documentclass[10pt]{article}
\usepackage[T1]{fontenc}
\usepackage{amsmath,amsfonts,amsthm,mathrsfs,amssymb}
\usepackage{cite}
\usepackage{multirow}
\usepackage{graphicx,float}
\usepackage{subfigure}
\usepackage{placeins}
\usepackage{color}
\usepackage{indentfirst}
\usepackage[bookmarks=true]{hyperref}
\definecolor{r1}{rgb}{1.000,0.000,0.000}
\numberwithin{equation}{section}
\newcommand{\be}{\begin{eqnarray}}
\newcommand{\ben}{\begin{eqnarray*}}
\newcommand{\en}{\end{eqnarray}}
\newcommand{\enn}{\end{eqnarray*}}

\newtheorem{theorem}{Theorem}[section]
\newtheorem{lemma}[theorem]{Lemma}

\newtheorem{remark}[theorem]{Remark}

\newtheorem{algorithm}{Algorithm}[section]

\definecolor{rot}{rgb}{1.000,0.000,0.000}

\usepackage{algorithm}   
\usepackage{algorithmic}  
\usepackage{subfigure}
\usepackage{lipsum}
\usepackage{amsmath}
\usepackage{amssymb}
\usepackage{bm}
\usepackage{pifont} 

\newcommand{\figref}[1]{Fig.~\ref{#1}}

\title{PFWNN: A deep learning method for solving forward and inverse problems of phase-field models}

\author{Gang Bao\thanks{School of Mathematical Sciences, Zhejiang University, Zhejiang 310058, China. Email: {\tt baog@zju.edu.cn}}\;,
Chang Ma\thanks{School of Mathematical Sciences, Zhejiang University, Hangzhou, Zhejiang 310058, China. Email:{\tt 12035014@zju.edu.cn}}\;,
Yuxuan Gong\thanks{School of Materials Science and Engineering, Zhejiang University, Hangzhou, Zhejiang 310058, China. Email: {\tt yxgong@zju.edu.cn}}}
 
\begin{document}

\maketitle

\begin{abstract}
    Phase-field models have been widely used to investigate the phase transformation 
    phenomena. 
    However, it is difficult to solve the problems numerically 
    due to their strong nonlinearities and higher-order terms.
    This work is devoted to solving forward and inverse problems of the phase-field 
    models by a novel deep learning framework named Phase-Field Weak-form 
    Neural Networks (PFWNN),
    which is based on the weak forms of the phase-field equations.
    In this framework, the weak solutions are parameterized as deep neural networks with periodic layers, 
    while the test function space is constructed by functions compactly supported in small regions. 
    The PFWNN can efficiently solve the phase-field equations characterizing the sharp transitions and identify the important parameters by employing the weak forms.    
    It also allows local training in small regions, which significantly reduce the computational cost. Moreover, it can guarantee the residual descending along the time marching direction, enhancing the convergence of the method. 
    Numerical examples are presented for several benchmark problems.
    The results validate the efficiency
    and accuracy of the PFWNN. This work also 
    sheds light on solving the forward and inverse problems 
    of general high-order time-dependent partial differential equations.
\end{abstract}


\section{Introduction}
\label{introduction}
The phase-field models have emerged as powerful mathematical tools
for studying the spatio-temporal evolutions along with certain physical properties
in diverse fields of science and engineering.
In particular, they have been widely used in materials science 
for describing a variety of microstructural evolution phenomena 
including solidification \cite{kim2004phase,hotzer2015large}, 
grain growth \cite{krill2002computer, chang2017effect}, precipitate coarsening \cite{zhu2004three,cheng2021phase} , 
dislocation dynamics \cite{chan2010plasticity,beyerlein2016understanding}, 
and crack propagation \cite{aranson2000continuum, karma2001phase}.
According to the models, the phase-field dynamics are governed by the free energy of the system, 
and the phase-field variables are characterized by distinct constant values exhibiting 
smooth transitions at the interfaces. 

In this work, we consider
 two typical phase-filed models, namely the Allen-Cahn equation,
\begin{equation}
    \frac{\partial \phi  (t, x)}{\partial t} = - M  
    \frac{\delta \mathcal{F} }{\delta \phi} \quad in \ \Omega_T =\Omega \times (0,T],
    \label{PF_AC}
\end{equation}
and the Cahn-Hilliard equation, 
\begin{equation}
    \frac{\partial \phi  (t, x)}{\partial t} = M  \Delta
    \frac{\delta \mathcal{F} }{\delta \phi} \quad in \ \Omega_T =\Omega \times (0,T],
    \label{PF_CH}
\end{equation}
\noindent
with the initial condition
\begin{equation*}
    \phi (x,0) = \psi (x),\quad x \in \overline{\Omega}.
\end{equation*}
\noindent
Here, 
$T>0$ is a finite time, and $M$ is the mobility coefficient.
The domain $\Omega$ is assumed to be d-dimensional ($d=1,2,3$) cube 
$ [0, L_1] \times \cdots \times [0, L_d]$.
Throughout, for simplicity, periodic boundary conditions are imposed for both equations.

The total free energy functional $\mathcal{F}$ takes the following form: 
\begin{equation}
    \begin{aligned}
    \mathcal{F}(\mathrm{\phi}) =\int_{\Omega}^{} \left\{F (\phi) + 
    \frac{\lambda^2}{2}|\nabla \mathrm{\phi}|^2 \right\} \,d \Omega, 
    \label{total}
    \end{aligned}
\end{equation}
where $\lambda$ is a small parameter characterizing the width of the
phase transition interface, and $F(\phi)$ is a nonlinear 
energy potential.
In this case, (\ref{PF_AC}) and (\ref{PF_CH}) can be rewritten as
\begin{align}
    \label{AC_0}\frac{\partial \phi }{ \partial t} &= M(\lambda^2 \Delta \phi- f(\phi)),\\    
    \label{CH}\frac{\partial \phi }{ \partial t} &= -M(\lambda^2 \Delta^2 \phi- \Delta f(\phi)),
\end{align}
respectively, where $f(\phi)=F'(\phi)$. Furthermore, it is assumed that $F(\phi)$ is taken as the double-well polynomial in this work,
\begin{equation}
    F(\phi) = \frac{1}{4} (\phi-1)^2(\phi+1)^2, \quad
    f(\phi) = \phi^3-\phi.
\end{equation}
In \eqref{AC_0}, $\phi$ is a non-conserved variable representing the disordered-ordered state during the dynamics.
In \eqref{CH}, $\phi$ represents the concentration of substance during the evolution process 
and is therefore a conserved variable.
Moreover, in certain applications of the materials science, these two equations 
can be interrelated or coupled as a system which describes the transition 
of various phase-fields simultaneously
\cite{WOS:A1994PL38400007,Barrett2001FINITE}.


During the past decades, a variety of traditional numerical methodologies 
have been employed to solve the forward problems of phase-field models. 
Spatial discretization techniques such as finite difference methods, 
finite element methods, and Fourier-spectral methods 
have been extensively utilized in addressing these problems \cite{chen1998applications, chen2002phase, shen2010numerical}.
More attention has been paid to designing an efficient 
discretization scheme in the temporal domain for the high-ordered 
systems with strong nonlinearities, such as convex splitting \cite{eyre1998unconditionally}, 
exponential time differencing \cite{fu2022energy} 
and scalar auxiliary variable \cite{shen2018scalar}.
But these methods have high requirements for temporal and spatial discrete step sizes.
In addition, some important parameters of 
the phase-field models are not often available in the actual problems. 
It is essential
to identify these parameters from the observed phase-field data.
Several theoretical and numerical results 
concerning the parameter identification in nonlinear parabolic 
equations can be found in \cite{cannon1980inverse, cao2006natural}. 
Recent studies have investigated the theoretical frameworks and numerical approaches for identification of the multiple parameters in the Cahn-Hilliard equation 
\cite{brunk2023uniqueness} and the phase-field model for tumor growth \cite{kahle2019bayesian}.
However, they only conducted numerical analysis for the case in one dimension.
Therefore, significant
challenges remain in solving both forward and inverse problems of phase-field models because
of the strong nonlinearities and the perturbation of 
Laplacian or biharmonic operators due to
the small parameters (i.e., $\lambda\ll 1$).
Therefore, innovative and efficient numerical methods are needed to address the challenges posed by the high requirements for temporal and spatial mesh sizes, as well as high-dimensional inverse problems.

The goal of this work is to present a novel deep learning approach for 
addressing these challenges in solving the phase-field models.
In recent years, deep learning methods have shown promising performance in solving 
the forward and inverse problems of partial differential equations (PDEs),
which provide useful complements and extensions for traditional numerical algorithms. 
One of the most popular deep learning methods is 
the physics-informed neural networks (PINNs) 
\cite{raissi2019physics} based on the strong form of PDEs, 
whose loss function consists of PDE residuals 
along with initial and boundary conditions. 
In addition, PINNs are capable of estimating unknown parameters 
by optimizing the match between the PDE solutions and the observed data.
Recently, PINNs have been further extended 
to solve the Allen-Cahn and Cahn-Hilliard equations by introducing
various sampling strategies in both spatial and temporal domains in \cite{mattey2022novel,wight2020solving}. 
However, in order to achieve satisfactory accuracy, a massive number 
of collocation points have to be used, 
resulting in high training costs. 
In contrast to strong-form methods, the weak-form methods formulate the loss function 
based on the weak/variational formulation of PDEs \cite{ew2018deep,li2019d3m, CiCP-29-1365}.
These methods offer advantages including requiring less smoothness 
of the solutions, fewer quadrature points, and the ability to 
facilitate local learning through domain decomposition
\cite{Kharazmi2019VariationalPN}.
More recently, the weak adversarial networks 
\cite{zang2020weak, bao2020numerical} have been introduced, 
which convert the problem into an operator norm minimization 
problem to solve the forward and inverse problems of PDEs.
The weak solutions and the test functions are both parameterized as 
deep neural networks, where the parameters are learned by an 
adversarial training strategy.
In \cite{zang2023particlewnn}, 
the test function space is further constructed 
by functions compactly 
supported in extremely small regions centered around particles.
Other construction strategies 
of the test function space along this direction can be found in 
\cite{khodayi2020varnet, Kharazmi2019VariationalPN,kharazmi2021hp}.

In this work, we propose
the Phase-Field Weak-form Neural Networks 
(PFWNN), a novel weak-form deep learning-based framework, which is designed for 
solving the forward and inverse problems of the phase-field models. 
The approach is based on the weak forms of the Allen-Cahn and Cahn-Hilliard equations.
Here, the phase-field variable is simulated 
by a deep neural network $\phi^{NN}$, 
where the Fourier features are embedded in order 
to satisfy the periodic boundary conditions.
The test functions are localized, which set to be a 
series of compactly supported functions.
In addition, to facilitate the training of the PFWNN, a simple modification for the loss function based on the weak form is employed, which guarantees the residual descending along the time marching direction.
For the inverse problems, we focus on the identification of the unknown free energy potential $f(\phi)$ by the observed spatio-temporal data. The function is represented by 
another neural network $f^{NN}$, 
where $\phi^{NN}$ is the input data.
Our primary contributions include:
We propose a novel weak-form framework for solving the forward and inverse problems the phase-field equations by using deep learning method.
Our approach is discretization-free, highly parallelizable, and more effective in capturing the solution of the phase-field models.
We demonstrate the convergence and efficiency of the proposed framework with several numerical examples and show its superiority in solving the phase-field equations with complicated solutions and inverse problems for identifying the energy functional.

The remainder of the paper is organized as follows. 
In Section \ref{PFWNN_net}, the PFWNN and algorithms are
 proposed for solving 
the forward and inverse problems of the phase-field models, 
which are based on 
the weak formulations of the phase-field models.
It includes the construction of the weak forms and the training details of the PFWNN.
In Section \ref{experiment}, the accuracy and efficiency of 
the proposed method are examined through several numerical experiments. 
Specifically, we presented some cases of forward and inverse problems for one and two  dimensional of the Allen-Cahn and Cahn-Hilliard equations.
Finally, this paper is concluded with general discussions in Section \ref{conclusions}.

\section{Weak Neural Networks for phase-field Models}
\label{PFWNN_net}
\begin{figure}[!htbp]
    \centering  
    \includegraphics[width=\textwidth]{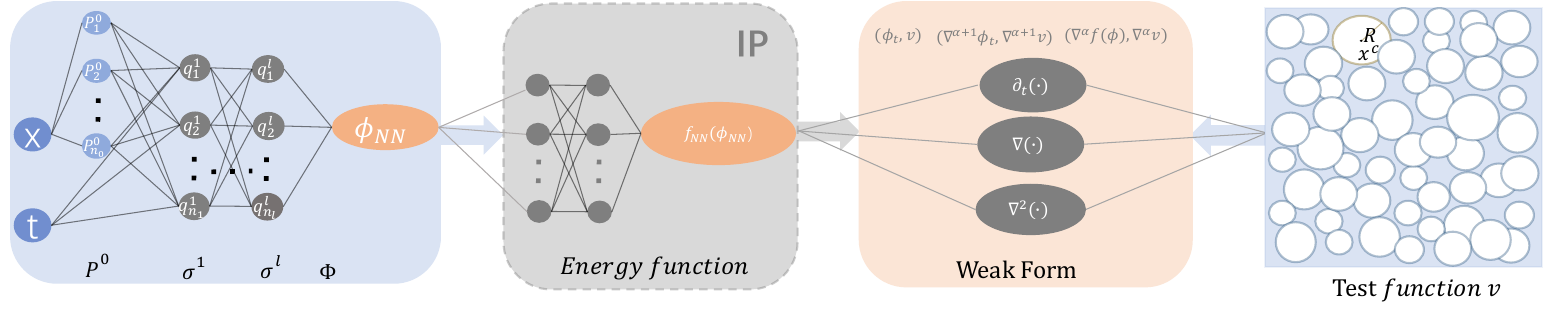}
    \caption{The sketch of the PFWNN. The dotted box should be added 
    when solving the inverse problems.} 
    \label{sketch}
\end{figure}

In this section, the framework of the PFWNN is introduced for solving the 
forward and inverse problems of the phase-field models. 
In the PFWNN, the loss function is based on the weak forms of the models, 
and the weak solutions are parameterized as deep neural networks 
with periodic structure. 
Motivated by our previous work \cite{zang2023particlewnn},
the test functions are selected as a series of functions 
compactly supported in small regions. In addition,
the free energy potential $f(\phi)$ 
is also represented as a deep neural network when solving the inverse problems by the PFWNN.
\figref{sketch} illustrates the idea of the PFWNN with a sketch.

\subsection{The Framework of the PFWNN}
Consider the general form of \eqref{AC_0} and \eqref{CH}:
\begin{equation}
    \frac{\partial \phi }{\partial t} =  M (- \Delta)^{\alpha} 
    (\lambda^2 \Delta \phi- f(\phi)), \quad \text{in}\ \Omega_T,
    \label{PF_AC_CH}
\end{equation}
where $\alpha =0$ and $1$ represent Allen-Cahn and Cahn-Hilliard equations, respectively.
After integration-by-parts and using the
periodicity of $\phi$ and $v$, the weak formulation of \eqref{PF_AC_CH} involves finding the solution
$\phi(x,\cdot) \in H^{\alpha+1}(\Omega;\mathbb{R})$, such that for all test
functions $v \in H_0^{\alpha+1}(\Omega;\mathbb{R})$, the following equation 
holds:

\begin{equation}\label{weak}
    (\phi_t, v) = - M \lambda^2 ( \nabla^{\alpha+1} \phi, \nabla^{\alpha+1} v)
    - M ( \nabla^\alpha f(\phi), \nabla^\alpha v), \quad t \in (0,T].
\end{equation}
\noindent
Here, $(\cdot, \cdot)$ is the standard $L^2$-inner product, 
i.e., $(\phi, v):=\int_\Omega \phi \cdot v dx $. 
The space $H^{\alpha+1}_0(\Omega;\mathbb{R} )$ denotes a Hilbert space of 
functions that themselves and $(\alpha+1)$-weak partial derivatives are $L^2$ 
integrable with vanishing trace on the boundary $\partial \Omega$.

It should be pointed out that as the differentiation order 
increases, more memory and computational cost 
for the automatic differentiation are inevitably required. 
This often results in a bottleneck when solving the equation 
with differential order higher than two \cite{wang2022respecting}. 
In order to save the computing resources, when $\alpha = 1$ , \eqref{weak} is reformulated as the following 
system of two coupled equations:
\begin{equation}
    \begin{aligned}
    &\phi_t = M\Delta \mu,  \quad \text{in}\ \Omega_T,\\
    &\mu = f - \lambda^2\Delta \phi,  \quad \text{in}\ \Omega_T.
    \end{aligned}
    \label{CH_form}
\end{equation}
The corresponding weak form of 
the problem (\ref{CH_form}) is to find 
$(\phi,\mu)\in H^1\times H^1$, such that for all 
$(v_1, v_2) \in H_0^1 \times H_0^1 $, the following equations hold:
\begin{equation}
    \begin{aligned}
    &(\phi_t, v_1) = - M(\nabla \mu, \nabla v_1),\quad t \in (0,T],\\
    &(\mu, v_2) = (f , v_2) + \lambda^2(\nabla \phi, \nabla v_2),\quad t \in (0,T].
    \label{CH_mu}
    \end{aligned}
\end{equation}

\subsection{Neural Networks Representing Weak Solutions}

The PFWNN approximates the phase-field variable $\phi$ with 
a neural network $\phi^{NN}$, 
which is usually comprised of $l$ hidden layers with $n_i$ neurons in each 
layer with a nonlinear and smooth activation function $\sigma(\cdot )$. 
In order to assure that $\phi^{NN}$ automatically satisfies the periodic boundary conditions, 
a periodic structure is imposed on the first layer of $\phi^{NN}$, 
which is called the periodic layer. 
This method is based on the following property about Composite function involving the periodic functions:

\begin{lemma}\label{lema_per}
    ($C^{\infty}$ periodic conditions \cite{dong2021method}) 
    Let $p(x)$ be a given smooth periodic function with period $L$ on the real axis, 
    i.e. $p(x +L) = p(x)$ for all $x \in (-\infty, \infty)$, and $f(x)$ denote an arbitrary 
    smooth function. Define $\phi(x) =f(p) =f(p(x))$. Then 
    $\phi^{(i)} (x) = \phi^{(i)} (x z+ L), \ i = 0, 1, 2,\cdots$.
\end{lemma}

By Lemma 1, in the one-dimensional (1D) case, we construct a Fourier 
feature with period $L$ and embed it into 
the first layer of $\phi^{NN}$ as follows:
\begin{equation}
    \begin{aligned}
    &p(x)=(1, sin(wx), cos(wx), \cdots, sin(mwx), cos(mwx)), ~~
    w=\frac{2\pi}{L} \\
    &P^{\left(0\right)}(x)=\sigma(p(x)),
    \end{aligned}
    \label{1D_per}
\end{equation}
\noindent
where $m \in \mathbb{N}^+$ is an adjustable training hyperparameter. 
Here, $P^{(0)}$ is the periodic layer with $n_0$ neurons with $n_0 = 2m+1$ in one dimension. 
Due to the nonlinearity of $\sigma(\cdot)$,
this periodic layer contains not only the frequencies 
$kw$ with $k = 0,1,\ldots, m$, but also higher ones with the same period $L$.
Similar to \eqref{1D_per}, in the two-dimensional (2D) case, 
the Fourier features are encoded as
\begin{equation}\label{2D_fourierfeature}
    p(x,y)=
        \begin{bmatrix}
        \sin\left(\omega_xx\right)\sin\left(\omega_yy\right),\ldots,
        \sin\left(m_x\omega_xx\right)\sin\left(m_y\omega_yy\right)\\
        \sin\left(\omega_xx\right)\cos\left(\omega_yy\right),\ldots,
        \sin\left(m_x\omega_xx\right)\cos\left(m_y\omega_yy\right)\\
        \cos\left(\omega_xx\right)\cos\left(\omega_yy\right),
        \ldots,\cos\left(m_x\omega_xx\right)\cos\left(m_y\omega_yy\right)\\
        \cos\left(\omega_xx\right)\sin\left(\omega_yy\right),\ldots,
        \cos\left(m_x\omega_xx\right)\sin\left(m_y\omega_yy\right)
    \end{bmatrix}
\end{equation}
where $ M_x, M_y\in\mathbb{N}^+$. Here, $L_x, L_y$ 
are the periods of \eqref{2D_fourierfeature} in the $x$ and $y$ directions, respectively,
and $w_x=2\pi/L_x, w_y=2\pi/L_y $.
Besides, the deep neural networks often fit target functions from low to high frequencies,\cite{xu2019frequency}
but the solution of the phase field equation is sharp at the interface. Introducing high frequency information in advance is beneficial for the learning of the neural network.

\par
Thus, the weak solution $\phi^{NN}$ in the PFWNN
takes the following form,
\begin{equation}
    \phi^{NN}(X;\theta) = \Phi \circ \sigma^{\left(l\right)} \circ 
    \sigma^{\left(l-1\right)}\circ\cdot \cdot \cdot \circ \sigma^{\left(1\right)}
    \circ P^{\left(0\right)},
\end{equation}
\noindent
where $X$ is the input data of the training network, $\theta$ is the training parameters
 of the network, $\sigma^{(i)}(i=1,2,\cdots,l$) and $\Phi$ are the nonlinear and 
linear mappings with weights and biases to be trained, respectively.
Note that the periodic layer $P^{(0)}$ only acts on the spatial data.

\subsection{Calculation of the Loss Functions}
According to \eqref{weak}, 
we denote $\mathcal{R}^t(\phi^{NN},v)$ as
the weak-form residual with a fixed time $t$:
\begin{equation}\label{residual_weak}
    \mathcal{R}^t(\phi^{NN};v)
    =(\phi^{NN}_t, v) + M \lambda^2 ( \nabla^{\alpha+1} \phi^{NN}, \nabla^{\alpha+1} v)
    + M ( \nabla^\alpha f(\phi^{NN}), \nabla^\alpha v).
\end{equation}
Specially for the form of  the Cahn-Hilliard equation \eqref{CH_mu}, we have 
\begin{equation}\label{residual_weak_chmu}
    \begin{aligned}
    &\mathcal{R}^t_1(\phi^{NN};v)
    = (\phi^{NN}_t, v) + M(\nabla \mu(\phi^{NN}), \nabla v),\\
    &\mathcal{R}^t_2(\phi^{NN};v)
    = (\mu(\phi^{NN}), v) - (f(\phi^{NN}) , v) + \lambda^2(\nabla \phi^{NN}, \nabla v).
    \end{aligned}
\end{equation}
where $  \mu(\phi^{NN})= f(\phi^{NN}) - \lambda^2\Delta \phi^{NN}$. For the sake of simplicity, we use the symbol $\mathcal{R}^t$ to represent $\mathcal{R}^t$ in \eqref{residual_weak} and $\mathcal{R}^t_1$, $\mathcal{R}^t_2$ \eqref{residual_weak_chmu} uniformly.

In order to enforce the weak solution neural network 
to focus on extremely local regions rather 
than integrate over the entire domain, 
the test functions are 
chosen to be compactly supported 
and defined in small spheres $\{B(\bm{x}^c, R)\subset \Omega\}$. 
Here, $\bm{x}^c$ is a particle in $\Omega$
and $R$ is the radius of the sphere. 
In this setting, we avoid the burden of interface treatment  as we only have a single loss function on a local region.
In this work, we choose the
compactly supported radial basis functions (CSRBFs) 
as test functions, which are defined in $\{B(\bm{x}^c, R)\subset \Omega\}$ as
\begin{equation} \label{csrbf}
    v(r)=\left\{\begin{array}{l}
    v^+(r), 
    \quad r<1\\
    0, \quad r \geq 1
    \end{array},\right. \quad
    r(\bm{x})=\frac{\|\bm{x}-\bm{x}^c \Vert }{R}.
\end{equation}

We randomly generate $N_p$ particles $\{\bm{x}_i^c\}_{i=1}^{N_p}$ 
and the corresponding $\left\{R_i\right\}_{i=1}^{N_p}$, 
and define CSRBFs $\{v_i\}_{i=1}^{N_p}$ in each small 
neighbourhood $\left\{ B(\bm{x}_i^c, R_i)\right\}_{i=1}^{N_p}$.
Then, we obtain the mean square error of the weak-form residuals with 
a fixed time $t$:
\begin{equation}\label{eq:SpatioLoss}
    \mathcal{L}_{\mathcal{R}^t}=\frac{1}{N_{p}} \sum_{i=1}^{N_{p}}\left|\mathcal{R}^t\left(\phi^{N N} ; 
        v_i\right)\right|^2,
\end{equation}
where
\begin{equation}\label{eq:residual_weak_BxR}
    \begin{split}
        \mathcal{R}^t\left(\phi^{NN} ; v_i\right) = & 
        (\phi^{NN}_t, v_i)_{B(\bm{x}_i^c, R_i)}
        + M \lambda^2 ( \nabla^{\alpha+1} \phi^{NN}, \nabla^{\alpha+1} v_i)_{B(\bm{x}_i^c, R_i)} \\
        & + M ( \nabla^\alpha f(\phi^{NN}), \nabla^\alpha v_i)_{B(\bm{x}_i^c, R_i)}.
    \end{split}
\end{equation}
To simplify the calculation, we convert the integral regions
$B(\bm{x}_i^c, R_i)$ into the unit sphere $B(\bm{0},1)$ by a simple 
coordinate transformation. By generating $N_{int}$ collocation points 
$\{s_j\}_{j=1}^{N_{int}}$ from $B(\bm{0},1)$
and denoting $\bm{x}^i_j :=\bm{s}_jR_i + \bm{x}^c_i$, we approximate
\eqref{eq:residual_weak_BxR} by
\begin{equation}\label{eq:residual_weak_approx}
    \begin{split}
        &\mathcal{R}^t\left(\phi^{NN} ; v_i\right) 
    \approx  R^d_i \mathcal{V} (d)\sum_{j = 1}^{N_{int}} 
    \Big[ 
    \phi^{NN}_t(t, \bm{x}^i_j)v_i(\bm{x}^i_j)\\
     & + \frac{M \lambda^2 }{R^{\alpha+1}} \nabla^{\alpha+1}\phi^{NN}(t, \bm{x}^i_j)\nabla^{\alpha+1}_r v_i(\bm{x}^i_j) 
    + \frac{ M}{R^\alpha} \nabla^\alpha f(\phi^{NN}(t, \bm{x}^i_j)) \nabla^\alpha_r v_i(\bm{x}^i_j)\Big],  
    \end{split}
\end{equation}
where $\mathcal{V}(d)$ is the volume of the d-dimensional unit sphere.     
It should be noted that $R_i$ can not be too large or small. 
In the PFWNN, we adopt the R-descending strategy to generate 
$\{R_i$\} for each particle $\bm{x}_i^c$ from a range $[R_{min}, R_{max}]$, 
where $R_{max}$ gradually decreases with iterations until a 
given lower bound is encountered.

To conform the temporal causal structure that is inherent from the evolution 
of physical systems,
we adopt the approach described 
in \cite{wang2022respecting}. 
Suppose that $0=t_1<t_2<\cdots<t_{N_T}=T$ discretizes the 
temporal domain,
we define a weighted residual as 
\begin{equation}
    \mathcal{L}_{\mathcal{R}}  = \frac{1}{N_T}\sum_{k = 1}^{N_T} 
     w_k \mathcal{L}_{\mathcal{R}^{t_k}},
\end{equation}
where
\begin{equation}
    w_k = exp \left( -\epsilon \sum_{j=1}^{k-1} \mathcal{L}_{\mathcal{R}^{t_j}} \right),\ k=2,3,\cdots,N_T.
\end{equation}
Here, $\epsilon$ is a causality parameter that controls the steepness of 
the weights $w_k$. 
The weights are inversely exponentially proportional to the magnitude of the 
cumulative residual loss from the previous time steps.
It can be seen that at the beginning of the training, $w_k$ is too small so that $\mathcal{L}_{\mathcal{R}^{t_k}}$ is barely involved in the minimization of $\mathcal{L}_\mathcal{R}$. Throughout the rest of the training, all the previous residuals $\big\{{\mathcal{L}_{\mathcal{R}^{t_j}}}\big\}_{j=1}^{k-1}$ decrease to some value. Therefore, $\mathcal{L}_{\mathcal{R}^{t_k}}$ will be considered in the optimization since $w_k$ becomes large enough. As a result, the weak-form residual loss can be properly minimized along the time direction. This weighted method not only accelerates the network training but also yields the trained models with better accuracy.

In addition, the weak solution neural network $\phi^{NN}$ need to satisfy the initial condition $\psi$.
We generate $N_{init}$ collocation points $\left\{x_j\right\}_{j=1}^{N_{init}}$ 
in $\Omega$, and denote $\mathcal{L}_{\mathcal{I}}$
as the data mismatch on the initial condition:
\begin{equation}
    \mathcal{L}_{\mathcal{I} }=\frac{1}{N_{init}} \sum_{j=1}^{N_{init}}\left|\phi^{N N}(0,\bm{x}_j)
    -\psi(0,\bm{x}_j)\right|^2.
\end{equation}
\noindent
For the forward problems (\textit{FP}), 
the loss function of the PFWNN is formulated  as
\begin{equation}
    \mathcal{L}_{FP} = \omega_\mathcal{R} \mathcal{L}_\mathcal{R}  
    +\omega_\mathcal{I}  \mathcal{L}_\mathcal{I},
\end{equation}
where $\omega_\mathcal{R}$, $\omega_\mathcal{I}$
are weight coefficients of the weak-form residual and the initial data mismath, respectively.
Finally, the weak solution $\phi$ and the energy potential $f(\phi)$ 
can be learned by minimizing $\mathcal{L}_{FP}$. It should be noted that $\mathcal{R}$ consists of two parts for \eqref{residual_weak_chmu}.

\begin{algorithm}[t!]
    \caption{The PFWNN Algorithm}
    \label{PFWNN}
    \begin{algorithmic}[1]
    \STATE \textbf{Input:} 
    {$N_{T}$, $N_p$, $N_{int}$, $N_{init}$, $L$, $m$, $\epsilon$, 
    $\omega_\mathcal{R}$, $\omega_\mathcal{I}$, $\omega_\mathcal{D}$ }, 
    learning rate $\tau$, training iterations $Iters$.
    IP: $ N_s, \bm{x}^{sen}_j, t^{sen}_j, \mathcal{D}(\bm{x}^{sen}_j,t^{sen}_j)$.
    \STATE \textbf{Initialize:} Network architecture $\phi^{NN}$, 
    IP: Parameter network architecture $f^{NN}$.\WHILE{iterations $<Iters$ }
    \STATE {Generate integral points $\{(\widetilde{\bm{x}}^i_j ,t_k)| \widetilde{\bm{x}}^i_j \in B({0},1)\}$,
    radius $\{R_i\}^{N_p}_{i=1} \sim \text{Unif}[R_{min}, R_{max}]$, 
    particles $\{{\bm{x}}^c_i\} \sim \text{LHS} \footnotemark{} \{\bm{x} \in \Omega|dist(\bm{x},\varGamma )
    \geqslant R_i \}$, then, $\{\bm{x}^i_j = \widetilde{\bm{x}}^i_j*R_i +{\bm{x}}^c_i | {\bm{x}}^i_j \in B(\bm{x}_i^c,R_i)\}$,   
    inital data points $\{\bm{x}_j\}$. 
    IP: measurement dataset $\left\{(\bm{x}^{sen}_j,t^{sen}_j),\mathcal{D}(\bm{x}^{sen}_j,t^{sen}_j)\right\}$.}
    \STATE {Calculate the loss and update network parameters with the Adam optimizer: 
    $\theta \leftarrow \theta - \tau \nabla_\theta \mathcal{L}_{FP}$ (or $\mathcal{L}_{IP}$).}
    \ENDWHILE
    \STATE \textbf{Output:} {The weak solution network $\phi^{NN}$, IP: the identified parameter $f^{NN}(\phi^{NN})$.}
    \end{algorithmic}
\end{algorithm}

\footnotetext{Randomly sampled point locations are generated using a space 
filling Latin Hypercube Sampling (LHS) strategy.}

\begin{remark}\label{rm:CSRBF}
    In fact, there exist a variety of CSRBFs,
    such as Bump function \cite{FRY2002143}, Wendland's
    function \cite{wendland1995piecewise} and Wu's function \cite{wu1995compactly}.
    In the numerical examples presented in this work, we consider the following  
    Wendland's type CSRBFs:
    \begin{equation}
    v^+(r)=\frac{(1-r)^{l+2}}{3}\left[\left(l^2+4 l+3\right) r^2+(3 l+6) r+3\right],
    \end{equation}
    where $l=\lfloor d/2\rfloor +3$ and $d$ is the dimension of the domain.
    \figref{testfun} shows the plots of the test functions centered at 
    $x^c=0$ and $\bm{x}^c=(0,0)$ in 1D and 2D, respectively.
    \begin{figure}[!htbp]
        \centering  
        \includegraphics[width=0.8\textwidth]{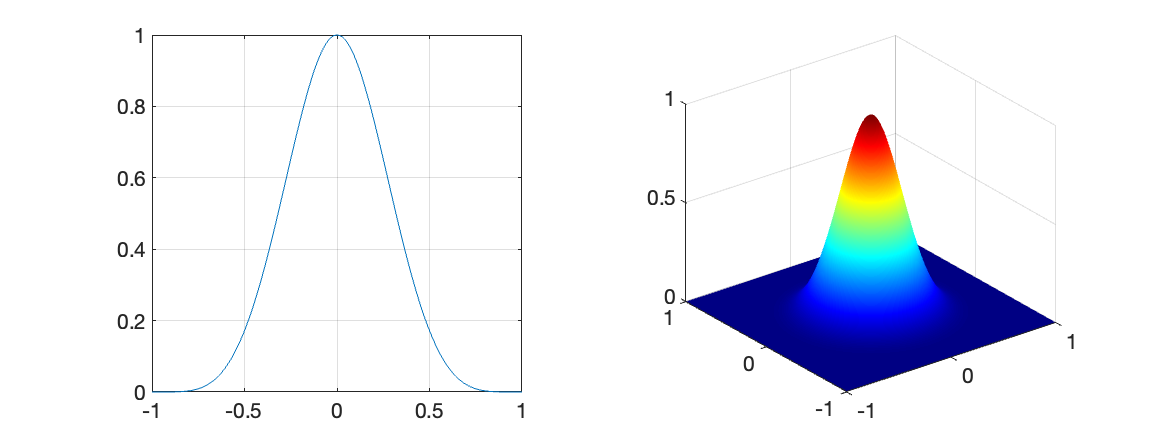}
        \caption{The plots of the test functions in 1D (left) and 2D (right).} 
        \label{testfun}
    \end{figure}
\end{remark}

\begin{remark}\label{rm:s2s}
    Recent work \cite{krishnapriyan2021characterizing,mattey2022novel,wight2020solving} demonstrates that it  
    may be more effective to formulate the forward problem as a sequence-to-sequence learning task, where the neural network predicts the solution for each temporal interval, instead of the entire temporal domain.
    Specifically, the entire temporal domain $[0, T ]$ is divided into sub-intervals
     $[0, \Delta t], [\Delta t, 2\Delta t], \cdots, [T -\Delta t, T ]$. Then, we train a network to learn the solution for each sub-interval, where the initial condition is obtained from the prediction of the previously trained network. At the end of training, the PFWNN is capable of predicting the target solutions over the entire spatio-temporal domain. They applied this idea to the PINN architecture, which we call it adaptive PINN(Ada-PINN).

\end{remark}

\subsection{Inverse Problems}

In the following, we study the identification of the energy functional $f(\phi)$. Assume that the mobility coefficient $M$ and the interface parameter $\lambda$ are known, and then expect that $f(\phi)$ can be identified uniquely from distributed measurements of  $\phi$ for the Allen-Cahn equations and Cahn-Hilliard equations.
From the Allen-Cahn equation \eqref{AC_0} and the Cahn-Hilliard equation \eqref{CH}, we obtain 
\begin{align}
     \label{AC_ff} M f(\phi) &= M\lambda^2 \Delta \phi - \phi_t,\\    
    \label{CH_ff} M\Delta f(\phi)&=\phi_t  +M\lambda^2 \Delta^2 \phi,
\end{align}
It is obvious that $f(\phi)$ can be uniquely determined in the Allen-Cahn equation \eqref{AC_ff} from observations of $\phi(\cdot,t)$ and $\phi_t(\cdot,t)$ on $\Omega$ for $t \in [0,T]$.
For the identification of $f(\phi)$ in the Cahn-Hilliard equation \eqref{CH_ff}, we can also obtain the uniqueness
 of the $f'(\phi)$ from the same observations, while $f(\phi)$ can be determined up to a constant shift \cite{brunk2023uniqueness}.

Therefore, we can identify the energy potential $f(\phi)$ from the measurement data $\phi(\cdot,t)$ and $\phi_t(\cdot,t)$. For the implementation of identifying $f(\phi)$ by the PFWNN,
we construct another neural network $f^{NN}$ where the weak solution neural network
$\phi^{NN}$ serves as the input, 
and the output is denoted by $f^{NN}(\phi^{NN})$. To obtain the measurement data $\phi_t(\cdot,t)$,
the size of time step is required to be small to approximate time derivatives.
Suppose we have $N_{s}$ sensors $\left\{(t^{sen}_j,\bm{x}^{sen}_j)\right\}^{N_s}_{i=1}$, 
and denote $\big\{\mathcal{D}(t^{sen}_j,\bm{x}^{sen}_j)\big\}$ as their noise measurements, 
the data mismatch $\mathcal{L}_{\mathcal{D}}$ on the sensors is defined as follows:
\begin{equation}
        \mathcal{L}_{\mathcal{D} }=\frac{1}{N_{s}} \sum_{j=1}^{N_{s}}\left|\phi^{N N}(t^{sen}_j,\bm{x}^{sen}_j)
        -\mathcal{D}(t^{sen}_j, \bm{x}^{sen}_j)\right|^2,
\end{equation}
Therefore, we obtain the loss function of 
the inverse problems (\textit{IP}),
\begin{equation}
        \mathcal{L}_{IP} = \omega_\mathcal{R}  \mathcal{L}_\mathcal{R} 
        +\omega_\mathcal{I}  \mathcal{L}_\mathcal{I}
        +\omega_\mathcal{D}  \mathcal{L}_\mathcal{D},
        \label{loss}
\end{equation}
where $\omega_\mathcal{D}$ is the weight coefficient of measurement data mismatch in the loss function.
In a nutshell, the detials for algorithm of the PFWNN are summarized 
in Algorithm \ref*{PFWNN}.

\begin{remark}\label{rm:CH_residual}
    If the form \eqref{CH_mu} is employed in solving the inverse problems of the Cahn-Hilliard equation 
    ($\alpha=1$), 
    an additional neural network  $\mu^{NN}$ can be constructed to represent $\mu$.
    The structure of the supplementary network is same as that of $\phi^{NN}$.
    In this way, we redefine 
    the weak-form residual as:
\begin{equation}\label{residual_weak_chmu_fnn}
    \begin{aligned}
    &\mathcal{R}(\phi^{NN},\mu^{NN};v)
    = (\phi^{NN}_t, v) + M(\nabla \mu^{NN}, \nabla v),\\
    &\mathcal{R}(\phi^{NN},\mu^{NN}, f^{NN};v)
    = (\mu^{NN}, v) - (f(\phi^{NN}) , v) - \lambda^2(\nabla \phi^{NN}, \nabla v).
    \end{aligned}
\end{equation}
    Here, three neural networks $\phi^{NN}$, $f^{NN}$ and $\mu^{NN}$ are trained simultaneously, and the structures of $\phi^{NN}$, $\mu^{NN}$ are the same.
\end{remark}

\section{Numerical Experiments}
\label{experiment}
Numerical results for solving the forward and inverse problems 
of the phase-field models are presented to validate the efficiency 
and accuracy of the PFWNN both in one and 
two dimensional cases.

\subsubsection*{Data generation}
To provide reference solutions for comparison with output of neural networks
and measurement data for inverse problems, 
we employ the semi-implicit Fourier-spectral method \cite{chen1998applications}.
By using the Fourier spectral approximation on spatio-domain, 
and representing the linear term implicitly and nonlinear term 
explicitly on temporal-domain, the discretization
forms of \eqref{AC_0} and \eqref{CH} are as follows:
\begin{equation}
    \begin{aligned}
    \frac{\left\{\phi(t+\Delta t)\right\}_{\bm{k}} -\left\{\phi(t)\right\}_{\bm{k}}}{\Delta t} 
    &=-k^{2\alpha} M \left\{\frac{d f}
    {d c}\right\}_{\bm{k}}-\lambda k^{2\alpha+2} M \left\{\phi(t+\Delta t)\right\}_{\bm{k}},\\
    \left\{\phi(t+\Delta t)\right\}_{\bm{k}} &=\frac{\left\{\phi(t)\right\}_{\bm{k}}-k^{2\alpha} M 
    \Delta t\left\{\frac{d f}{d c}\right\}_{\bm{k}}}{1+k^{2\alpha+2} M \Delta t },
    \end{aligned}
    \label{generate}
\end{equation}
where $\bm{k} = (k_1, \cdots, k_d)$ is a vector in the Fourier space, 
$k=|\bm{k}|$ is the magnitude of $\bm{k}$, $\{\cdot \}_{\bm{k}}$ 
represents the Fourier transform of the function.

\subsubsection*{Experimental setups}
We choose the ResNet framework with activation function 
\textit{Tanh}, 
and the Adam optimizer for updating the parameters of networks.
By default, the weak solution networks $\phi^{NN}$ and $\mu^{NN}$ comprise 3 hidden layers and 50 nodes per layer for 1D cases, while we employ 4 hidden layers with 100 nodes per layer for 2D cases.
And the network $f^{NN}$ comprises 3 hidden layers with 20 nodes per layer.
The parameters of the neural networks are randomly initialized 
by Kaiming Initialization. 
We set the initial learning rate $\tau = 0.001$ 
and apply the StepLR scheduler for adjusting the learning rate 
with the decay multiplication factor $\gamma  = 1 - 1/{Iters} $ for each iteration, and $Iters$ is the maximum number of iterations.
The period layer is of period $L = 2$ 
in each dimension, and the hyperparameter $m$ is set to be 5.
The radius $R_i$ is randomly sampled from $[R_{min}, R_{max}]$ for each region with $R_{max}$ decreasing linearly for each iteration, where $R_{min} = 10^{-6}, R_{max} = 10^{-4}$ for 1D cases and $R_{min} = 10^{-4}, R_{max} = 10^{-2}$ for 2D cases.
The ratio of the weights $\omega_\mathcal{R}, \omega_\mathcal{I}$ in the loss function is set to be $2:1$.
For the forward problems of the Allen-Cahn equation and the Chan-Hilliard equation, we split up 5 and 10 sub-intervals with 20000 iteration steps for each sub-interval.
For inverse problems, we estimate the energy functional $f$on the entire spatio-temporal domain with 50000 iteration steps and add $0.05\%$ noise to the measurement data. Besides, we use the $L^2$ relative error $\| \phi^{NN}-\phi \|^2_2 / \|\phi \|^2_2$,  
$\| f^{NN}(\phi^{NN})-f(\phi) \|^2 / \|f(\phi) \|^2$
and the maximum absolute error $max |\phi^{NN} - \phi|$
as evaluation metrics 
and run each example with 5 random 
seeds to obtain the mean value of the errors. 
The numerical simulations in the 1D cases are conducted using the device 
with NVIDIA GeForce RTX 2080 Ti GPU, 
while all other experiments in the 2D cases are conducted 
using NVIDIA A800 SXM4 80GB GPU.

\subsection{Numerical Results for the Allen-Cahn Equations}
\subsubsection{The 1D case }

 \begin{figure*}[htb]
    \begin{minipage}{0.8\linewidth}
    \subfigure[]{\includegraphics[width=0.32\linewidth]{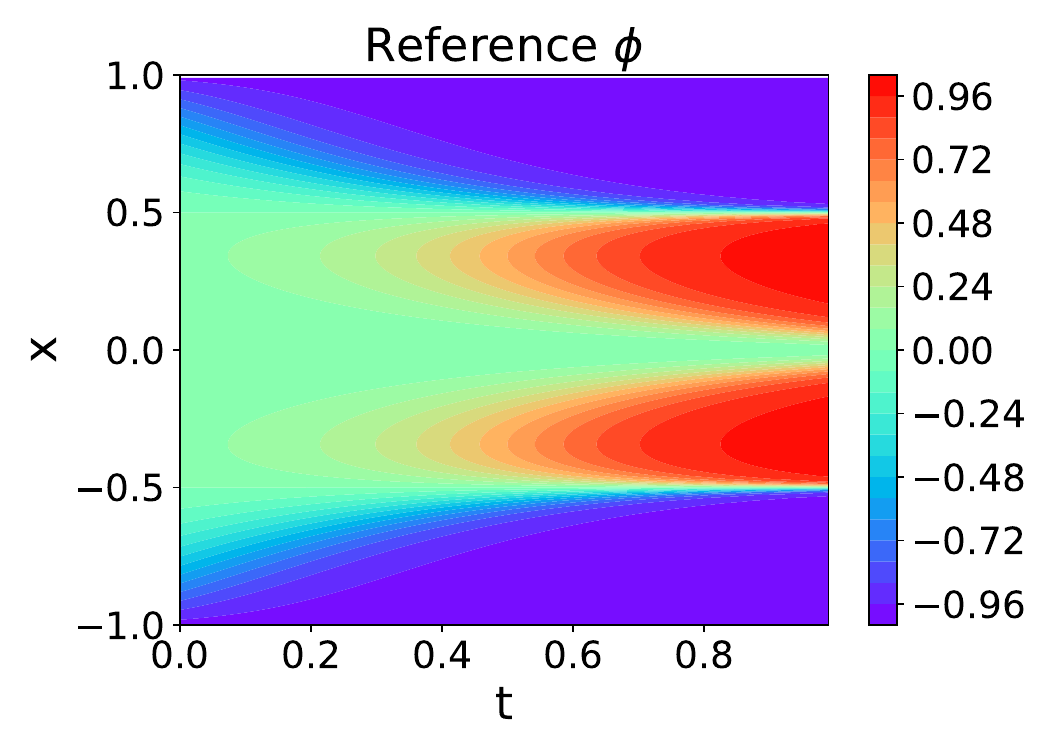}}
    \subfigure[]{\includegraphics[width=0.32\linewidth]{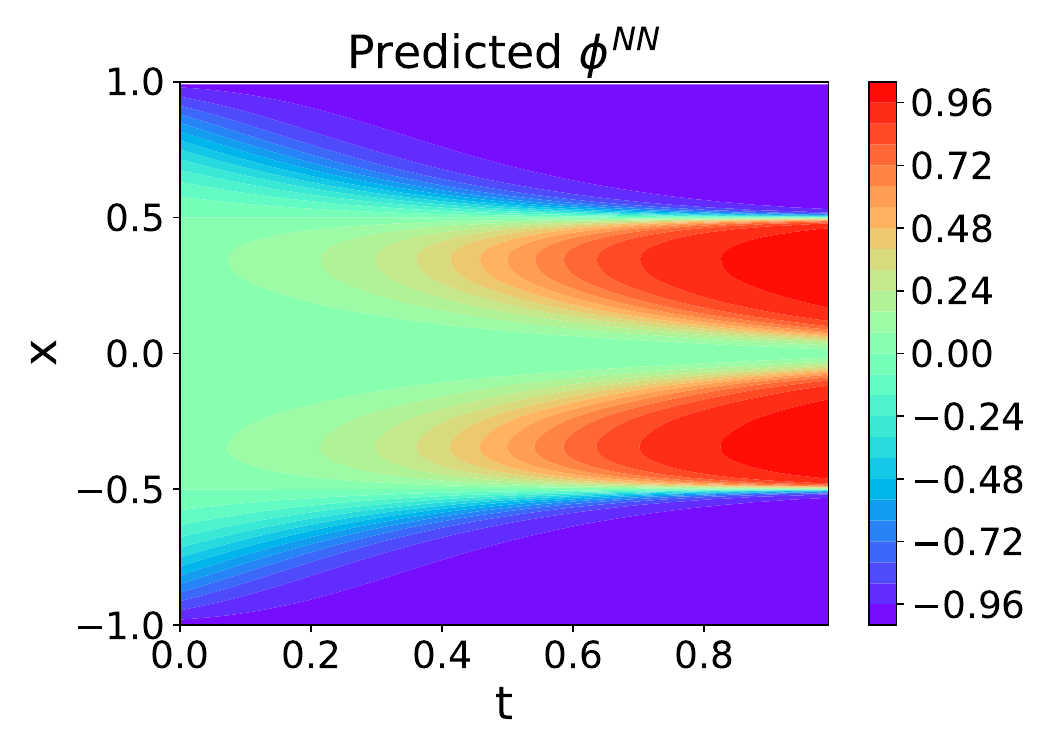}}
    \subfigure[]{\includegraphics[width=0.32\linewidth]{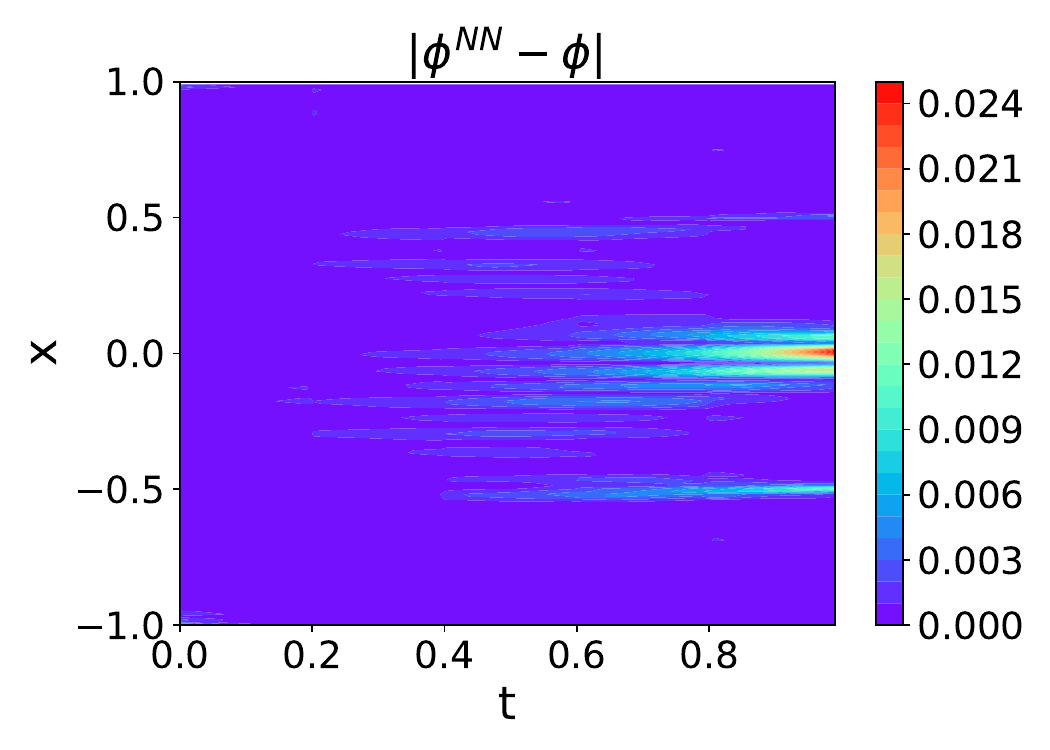}}
    \centering
    \end{minipage}
    \begin{minipage}{0.8\linewidth}
    \subfigure[]{\includegraphics[width=0.32\linewidth]{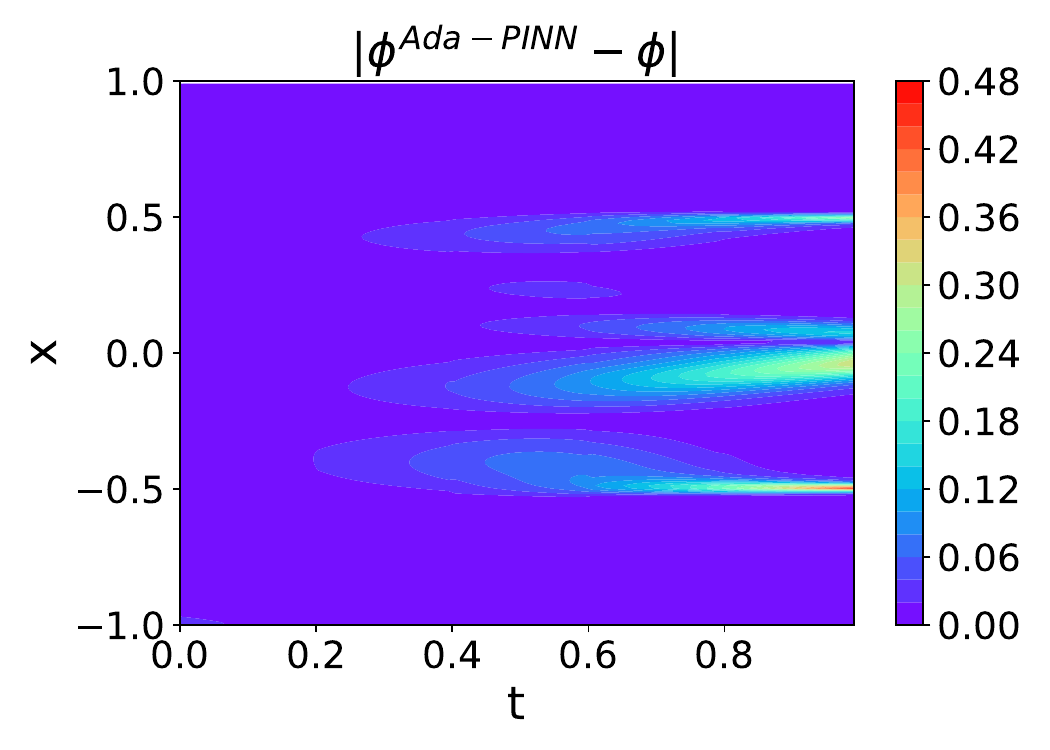}}
    \subfigure[]{\includegraphics[width=0.32\linewidth]{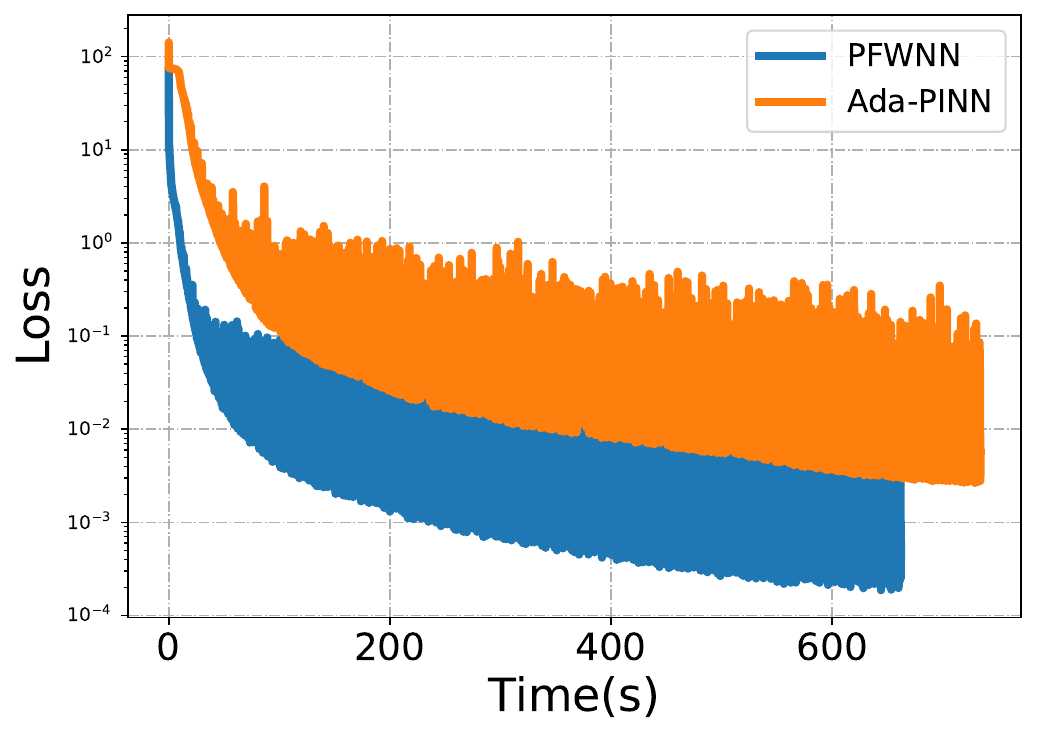}}
    \subfigure[]{\includegraphics[width=0.32\linewidth]{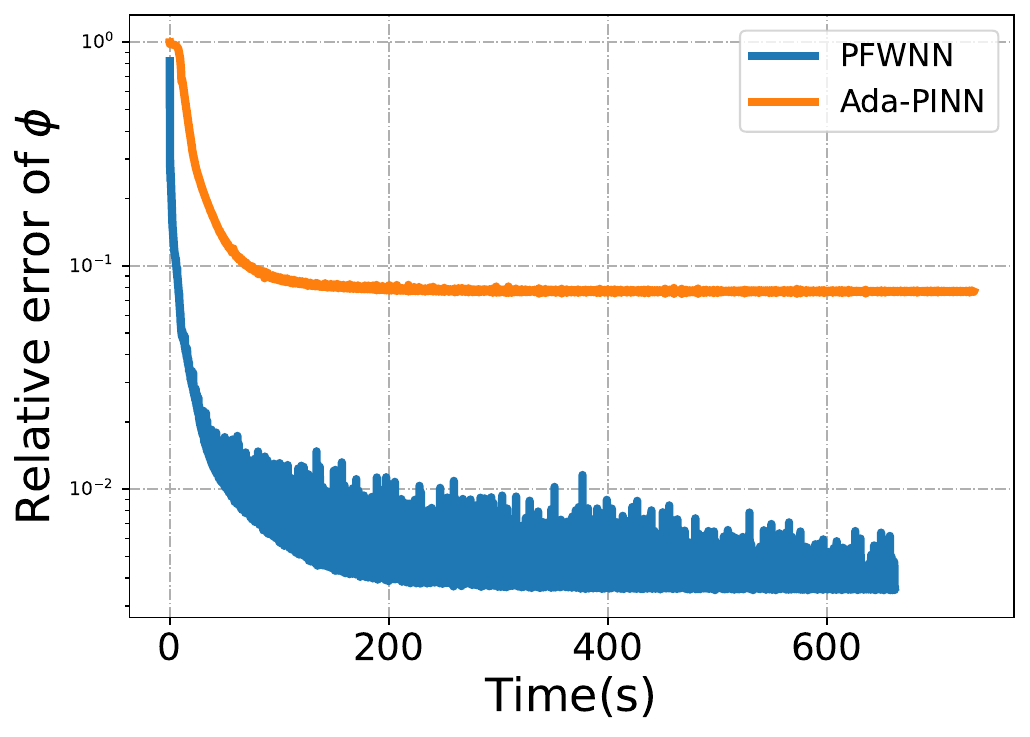}}
    \centering
    \end{minipage}
    \begin{minipage}{0.8\linewidth}
    \subfigure[]{\includegraphics[width=0.7\linewidth]{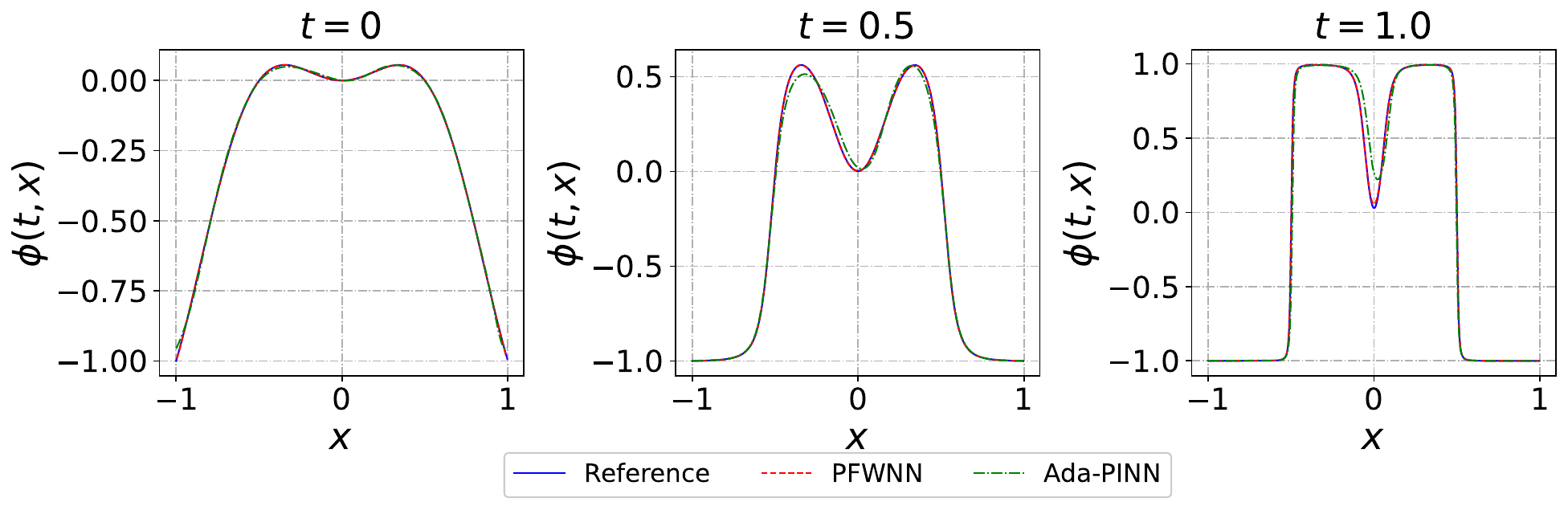}}
    \centering
    \end{minipage}
    \centering
    \caption{Results for 1D Allen-Cahn system. Spatio-temporal solutions:
            (a) The reference solution $\phi$, (b) The predicted solution 
            $\phi^{NN}$,
            (c) The point-wise error of the PFWNN.
            (d) The point-wise error of Ada-PINN.
            (e) Loss of the PFWNN and Ada-PINN vs. computation times for the last sub-interval.
            (f) Relative errors of solution $\phi^{NN}$ for the PFWNN and Ada-PINN vs. computation times for the last sub-interval.
            (g) The three plots are the predicted solutions of the PFWNN and Ada-PINN vs. reference solutions at different timestamps.\label{ac1dfp}} 
\end{figure*}
In this section, 1D time varying Allen-Cahn equation \eqref{AC_0} is considered. 
The values of the parameters are considered as
$M = 5, \lambda^2 = 2 \times10^{-5}$. 
The model simulates the phase separation phenomenon.
The final expression and initial conditions 
of the equation are as follows,
\begin{equation}\label{eq_ac1d_0}
    \begin{aligned}
    & \phi_t-\left(0.0001 \phi_{x x}-5\left(\phi^3-\phi\right)\right)=0, \ x \in[-1,1], \ t \in (0,1], \\
    & \phi(0,x)=x^2 \cos (\pi x), 
    \end{aligned}
\end{equation}

To evaluate integrals, we randomly sample $N_p = 50$ particles in the domain,
$N_{int}=5$ integration points in $B(0,1)$, $N_{T}=50$ temporal
collocation points and $N_{init}=100$ initial collocation points.
The causality parameter $\epsilon$ is set to be 0.1. 
The reference data is generated using \eqref{generate} 
with a spectral Fourier discretization 
of 512 modes and time-step size of $0.005$.
For comparison, we generate 20000 collocation points for the Ada-PINN to evaluate the strong residuals and keep other settings consistent with the PFWNN.
 
In \figref{ac1dfp}, we present the diagrams of the reference solution, 
the predicted solution and the pointwise error for the PFWNN and the Ada-PINN method.
The reference and the predicted solutions at three time points 
$t = 0, \ 0.5,\ 1.0$ are presented to further show the comparison. 
The resulting relative $L^2$ error of the PFWNN is $3.53e{-3}$, and the maximum absolute error of the PFWNN is $2.57e{-2}$. 
As a comparison, 
the resulting relative $L^2$ error and the maximum absolute error for the Ada-PINN are $7.48e{-2}$ and $4.92e{-1}$, respectively. 
The Ada-PINN most effectively captured solutions at the initial temporal domain, yet its accuracy decreases in subsequent intervals.
Specifically, the network does not perfectly learn the sharp curve as time approaches 1.0.
Even with a larger number of collocation points, the Ada-PINN does not achieve sufficient accuracy.
It indicates that the PFWNN achieves higher accuracy and faster convergence than the Ada-PINN, and shows that the deep learning method can relatively accurately predict the solution on the entire domain. 

Next, we further investigate the capability of the PFWNN for solving different problems. 
Note that the transition interface of the solutions is less sharp, 
which makes it easier to solve numerically when the parameter $\lambda$ increases or $M$ decreases. 
In particular, we change the initial condition and the parameters for the Allen-Cahn equation, i.e.,
\begin{equation}
    \begin{aligned}
    & \phi_t-\left(0.0001 \phi_{x x}-4\left(\phi^3-\phi\right)\right)=0, \ x \in[-1,1], \ t \in (0,1], \\
    & \phi(0,x)=x^2 \sin (2\pi x), 
    \end{aligned}
\end{equation}
\begin{figure*}[htb]
    \begin{minipage}{0.8\linewidth}
    \subfigure[]{\includegraphics[width=0.32\linewidth]{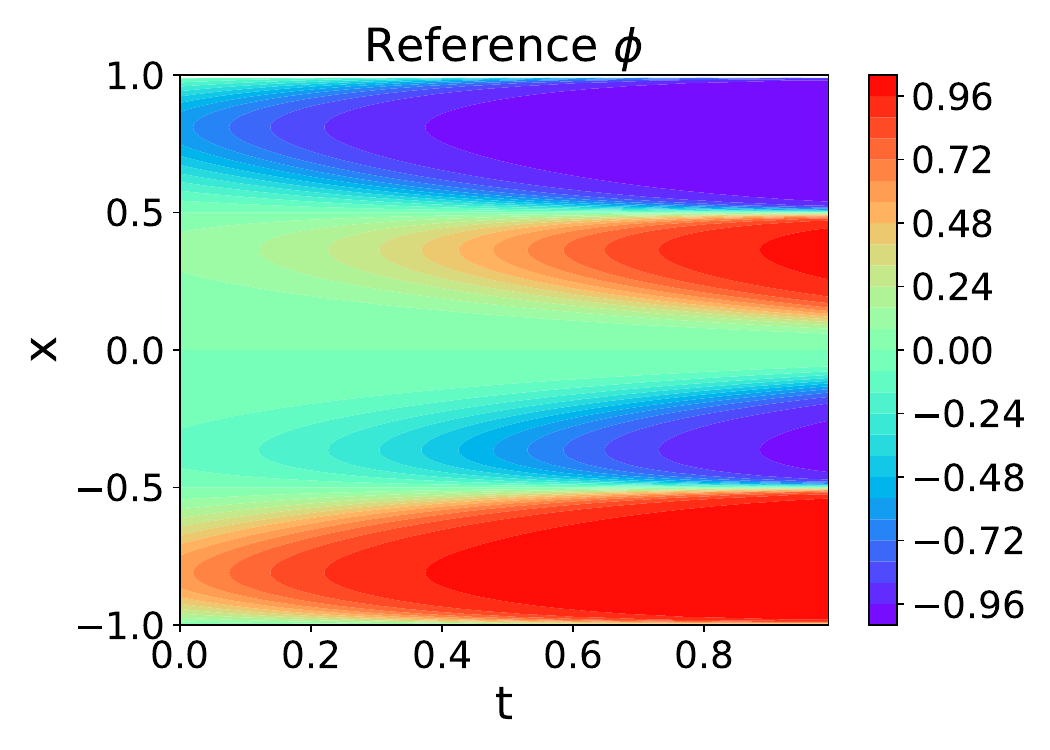}}
    \subfigure[]{\includegraphics[width=0.32\linewidth]{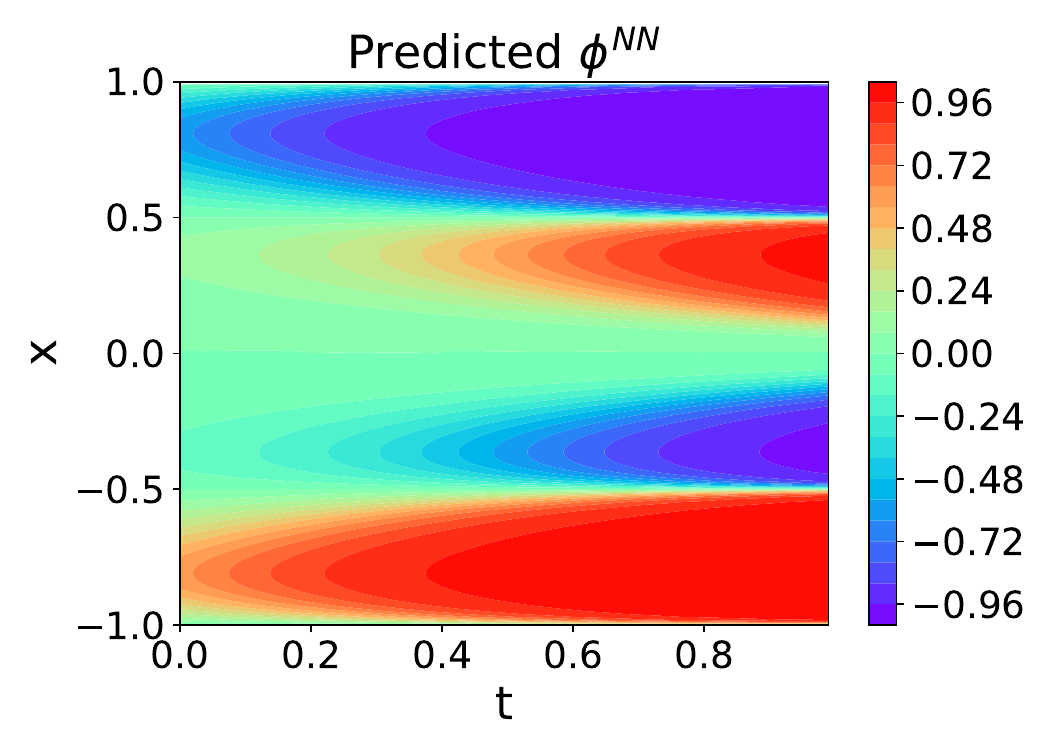}}
    \subfigure[]{\includegraphics[width=0.32\linewidth]{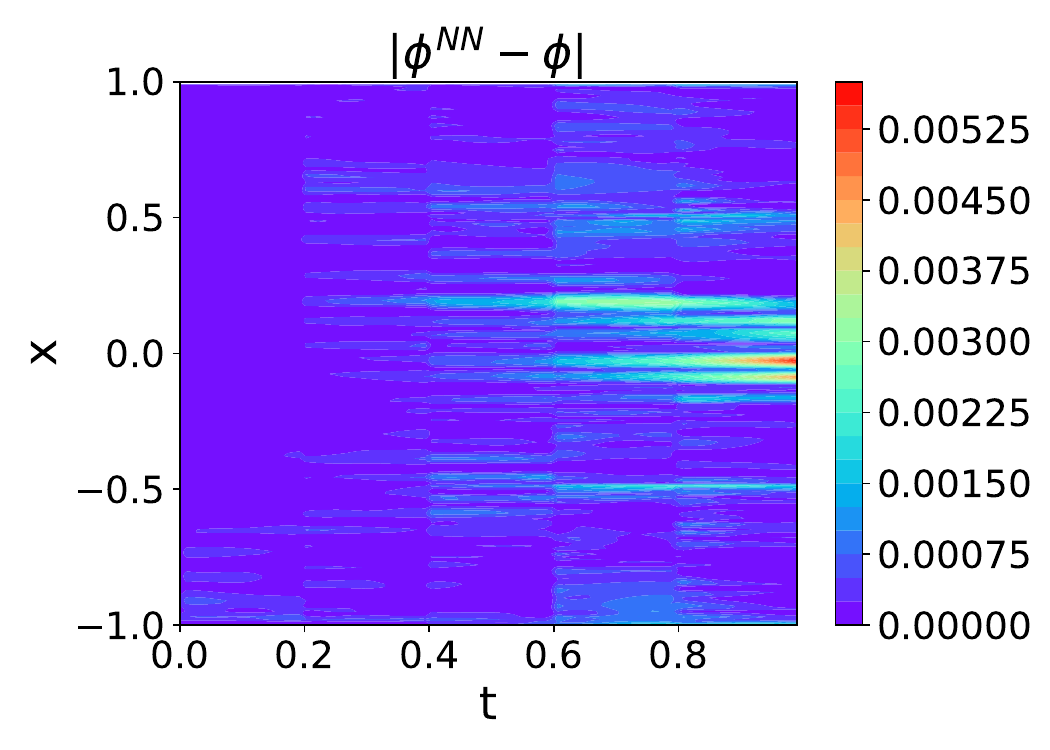}}
    \centering
    \end{minipage}
    \begin{minipage}{0.8\linewidth}
    \subfigure[]{\includegraphics[width=0.7\linewidth]{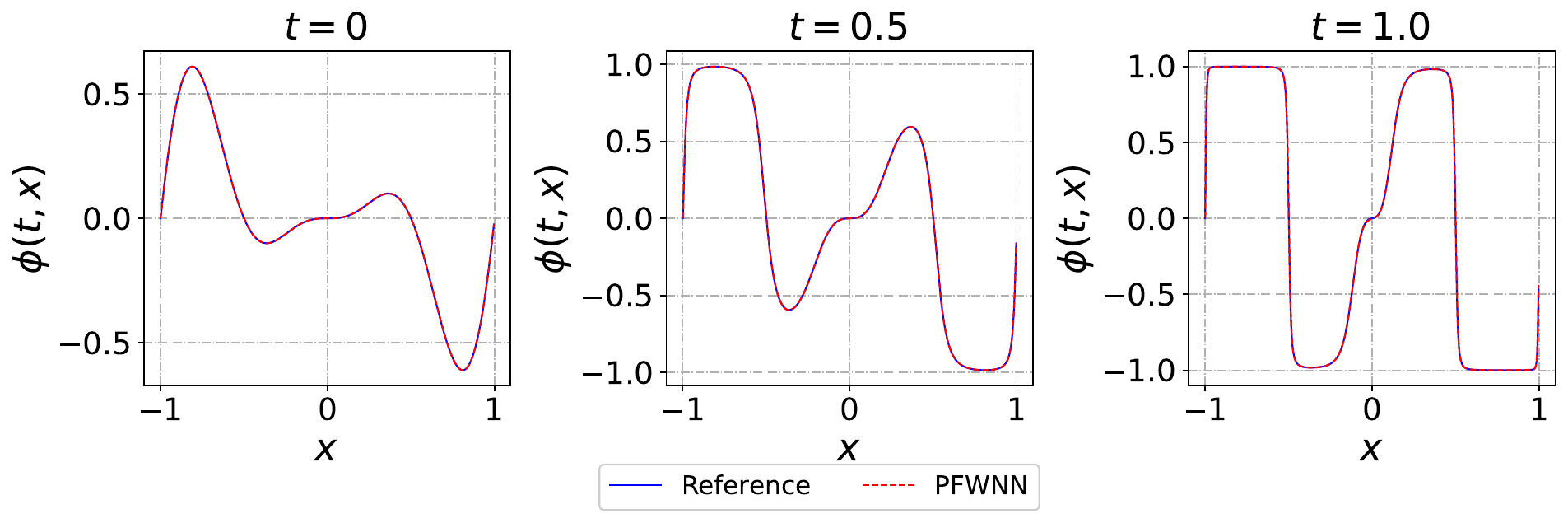}}
    \centering
    \end{minipage}
    \centering
    \caption{ Results for 1D Allen-Cahn system. Spatio-temporal solutions:
            (a) The reference solution $\phi$, (b) The predicted solution 
            $\phi^{NN}$,
            (c) The point-wise error of the PFWNN.
            (d) The three plots are the predicted solutions vs. reference solutions at different timestamps.\label{ac1dfp_1}} 
\end{figure*}
The neural network architecture and the parameter settings remain the same as 
that of the previous example.
In \figref{ac1dfp_1}, we observe that the result the previous example is lower than that of the example for the larger value of $M = 4$ and the smaller value of $\lambda^2 = 2.5\times10^5$. Therefore, we achieve higher accuracy with the relative $L^2$ error $1.07e{-3}$ and the maximum absolute error $5.65e{-3}$.

\begin{figure*}[!ht]
    \begin{minipage}{0.35\linewidth}
    \subfigure[]{\includegraphics[width=0.9\linewidth]{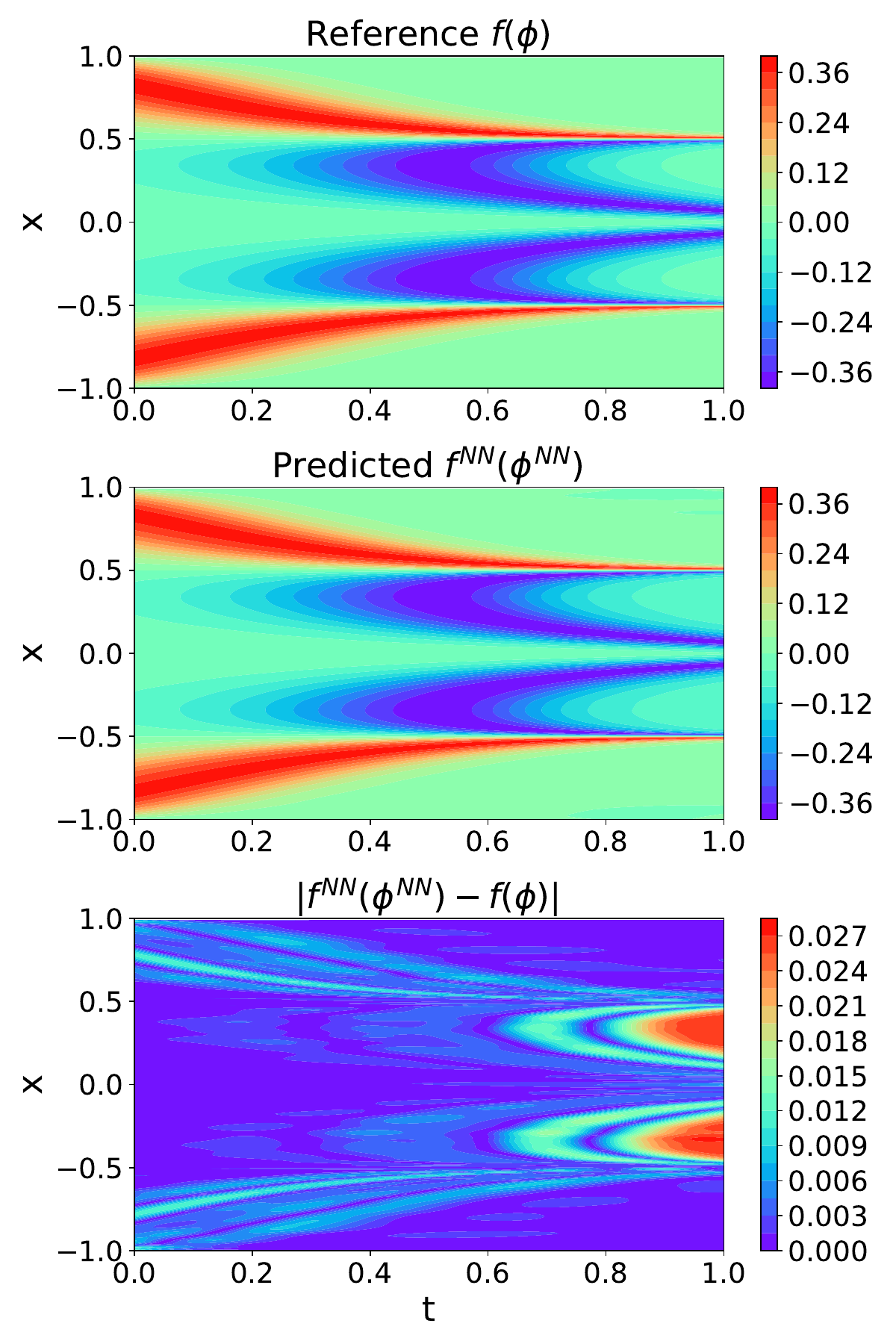}}
    \centering
    \end{minipage}
    \begin{minipage}{0.45\linewidth}
    \subfigure[]{\includegraphics[width=0.65\linewidth]{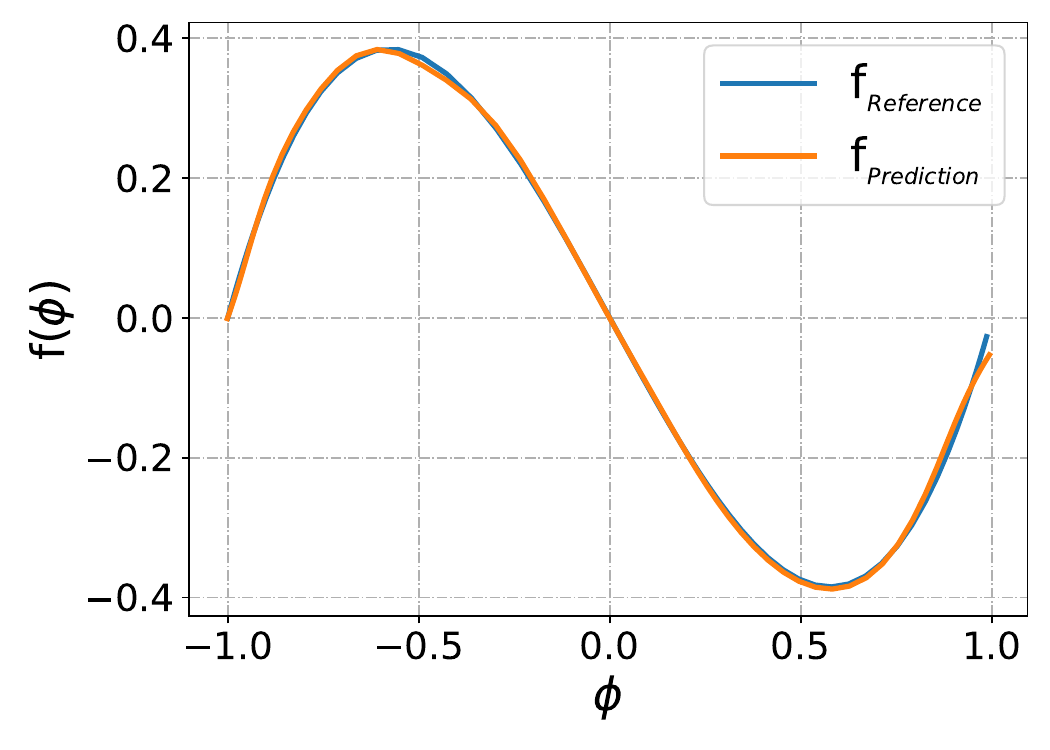}}
    \subfigure[]{\includegraphics[width=0.6\linewidth]{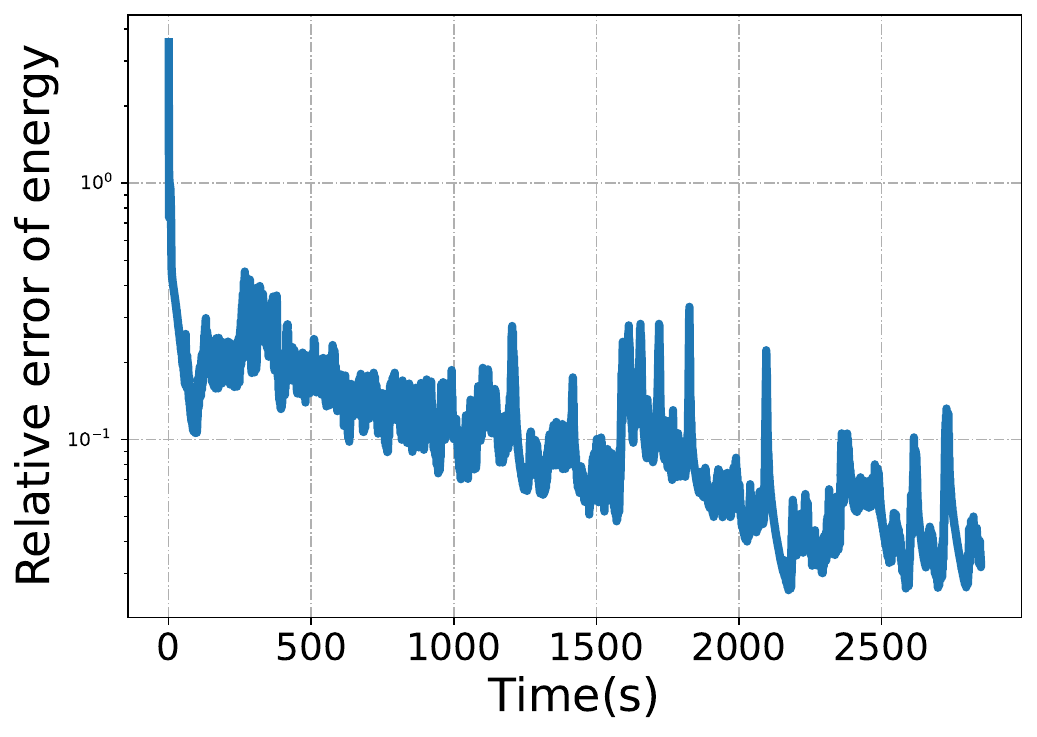}}
    \centering
    \end{minipage}
    \centering
    \caption{Results for the inverse problem of 1D Allen-Cahn system. 
            (a) Spatio-temporal energy functional:
            \textit{Top}: The reference energy functional $f(\phi)$, \textit{Middle}: The predicted solution 
            $f^{NN}(\phi^{NN})$,
            \textit{Bottom}: The point-wise error.
            (b) The predicted energy functional $f^{NN}$ vs. reference solutions $f$.
            (c) Relative errors for solution $f^{NN}$ vs. computation times.
             \label{ac1dip_1d}} 
\end{figure*}
For the inverse problem, we estimate the energy functional $f$ by the PFWNN on the entire domain. The other settings are consistent with the corresponding forward problem.
We obtain measurement data $\phi(\bm{x}, t_p)$ from the \eqref{eq_ac1d_0}.
The measurement data is input into the neural network $\phi^{NN}$, which serves as the input of the network $f^{NN}$.

Numerical results for this example are presented in \figref{ac1dip_1d}. We observe that the PFWNN can correctly identify the 
unknown parameters $f(\phi)$ with high accuracy even 
when the training data is corrupted with noise. 
The landscape 
The predicted energy function $f^{NN}$ is quite close to 
the reference energy function.
Specifically, the resulting relative $L^2$ error in estimating $f(\phi)$ is $9.56e{-4}$.
Consequently, the PFWNN produces stable and accurate reconstructions of the energy functional $f(\phi)$ only on the corresponding range of applicable data.

\subsubsection{The 2D case}
We study the evolution of the 2D Allen-Cahn equation which describes a single non-conserved phase particle.
Here, we choose the 2D cube $\Omega:=[-1,1]^2$ with 
the parameters $ M= 10 , \lambda^2 = 0.05^2$, and the initial profile for $\phi$ is given as
\begin{equation}\label{eq_ac2d_0}
    \begin{aligned}
    & \phi_t- 10\left( 0.05^2 \varDelta \phi- \left(\phi^3-\phi\right)\right)=0, \quad \bm{x} \in[0,1]^2, \quad t \in (0,T],\\
    & \phi( t=0, x, y)=\tanh \left(\frac{0.7-\sqrt{x^2+y^2}}{\kappa \lambda}\right),
\end{aligned}
\end{equation}
where $T = 3$, $\kappa = 2$.
During the training process, we set 
particles $N_p = 100$,
integration points $N_{int}=15$,
temporal collocation points $N_{T}=25$, 
iterations $Iters = 50000$, 
$N_{init}=500$ initial collocation points and 
causality parameter $\epsilon =1e-3$.
The reference data is generated by \eqref{generate} with a spectral Fourier discretization 
with $128 \times 128$  modes and time-step size of $10^{-2}$.

For identifying $f$ of \eqref{eq_ac2d_0}, 
We obtain measurement data $\phi(\bm{x}, t_p)$ from the reference solution where $\bm{x} \in \Omega ,\ t_p = 0.01p,\  p=0, \ 1,
\cdots, \ 10$. 
The noisy data is input into the neural network to calculate the data mismatch. In what follows, the resulting dataset corresponding to 
the predicted solution is used for model training, while the remaining 
data serves as the validation data.

\begin{figure*}[htb]
    \centering
    \begin{minipage}{0.6\linewidth}
    \subfigure[]{\includegraphics[width=\linewidth]{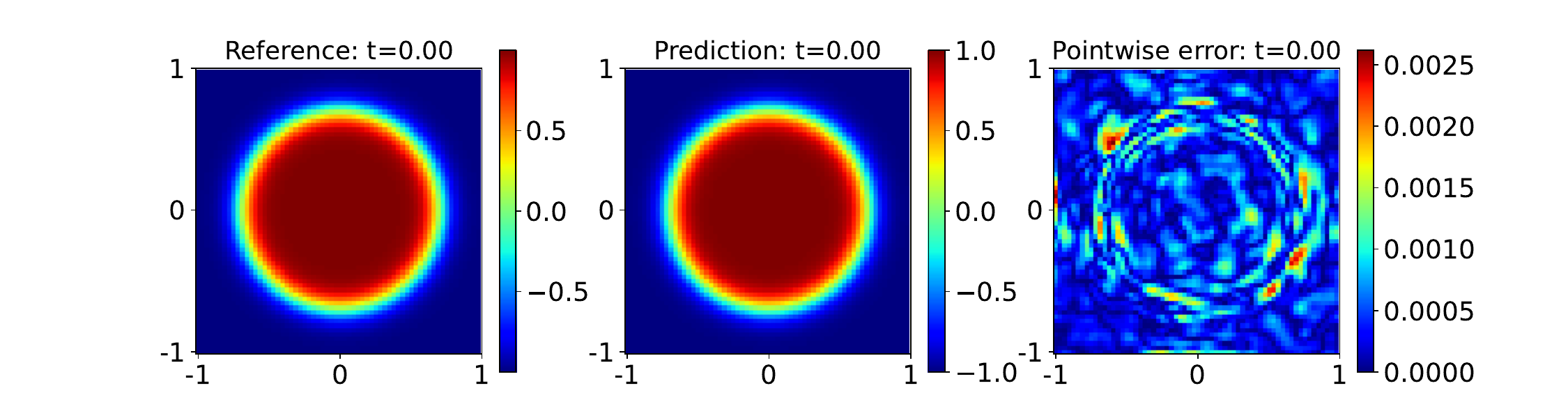}}
    \subfigure[]{\includegraphics[width=\linewidth]{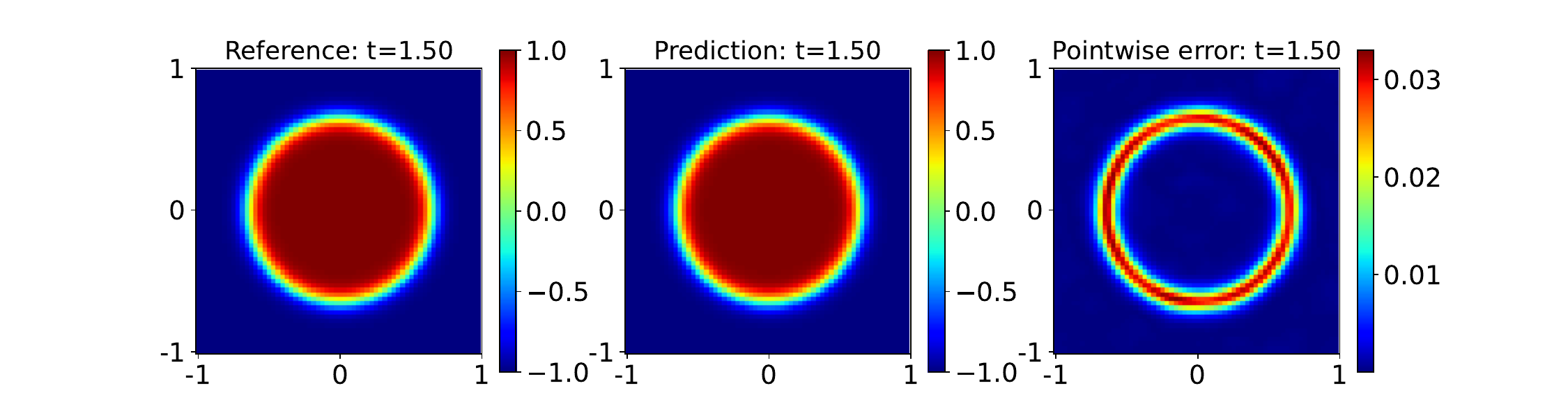}}
    \subfigure[]{\includegraphics[width=\linewidth]{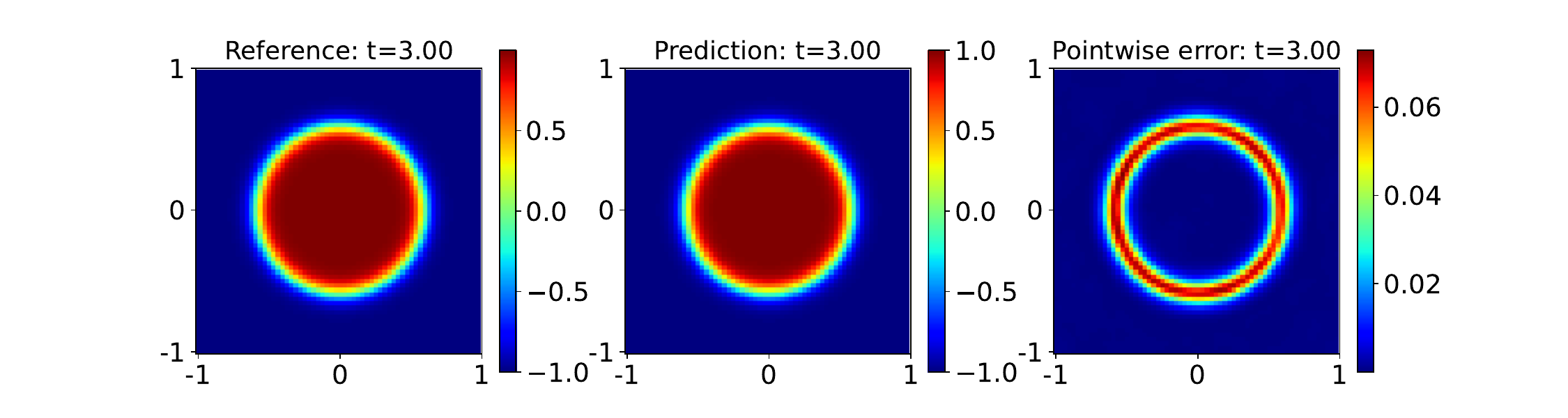}}
    \end{minipage}
    \caption{
    (a)-(c) The reference solutions and the predicted solutions 
    at different time snapshots.
    \label{ac2dip_u} }
\end{figure*}

A visual comparison between the reference solution and the predicted solution is presented in \figref{ac2dip_u}) at three time points $t=0, \ 1.5, \ 3.0$, where the phase particle shrink. It can be seen that the error is already very small. In other words, the PFWNN provides an accurate approximation for this problem. The resulting relative $L^2$ error and and the maximum absolute error of the PFWNN is $4.87e{-3}$ and $4.54e{-2}$, respectively. 

\begin{figure*}[htb]
    \centering
    \begin{minipage}{0.60\linewidth}
    \centering
    \subfigure[]{\includegraphics[width=\linewidth]{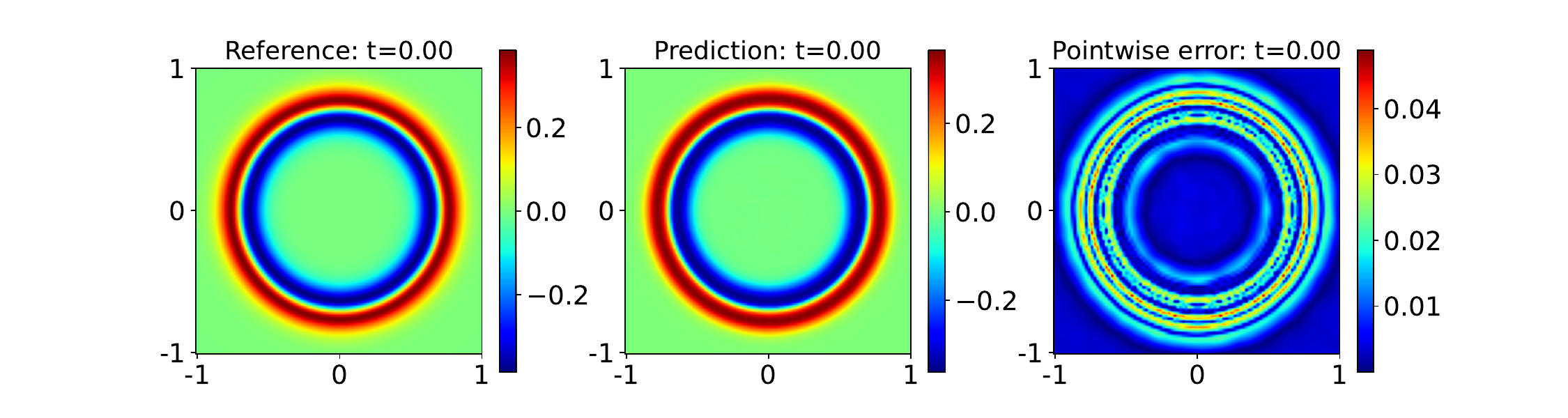}}
    \newline
    \subfigure[]{\includegraphics[width=\linewidth]{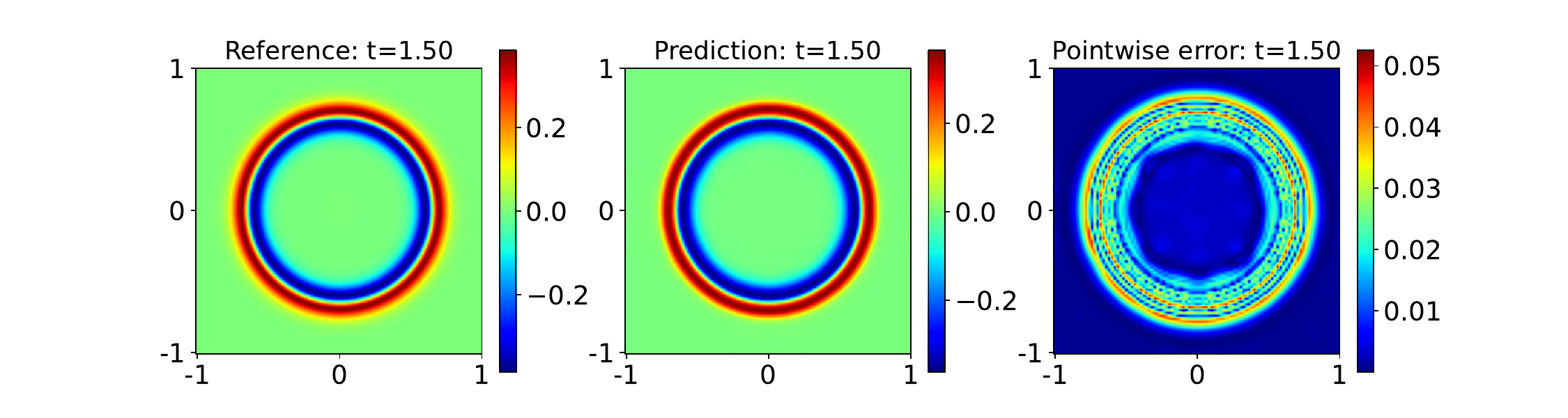}}
    \end{minipage}
    \begin{minipage}{0.35\linewidth}
    \subfigure[]{\includegraphics[width=0.65\linewidth]{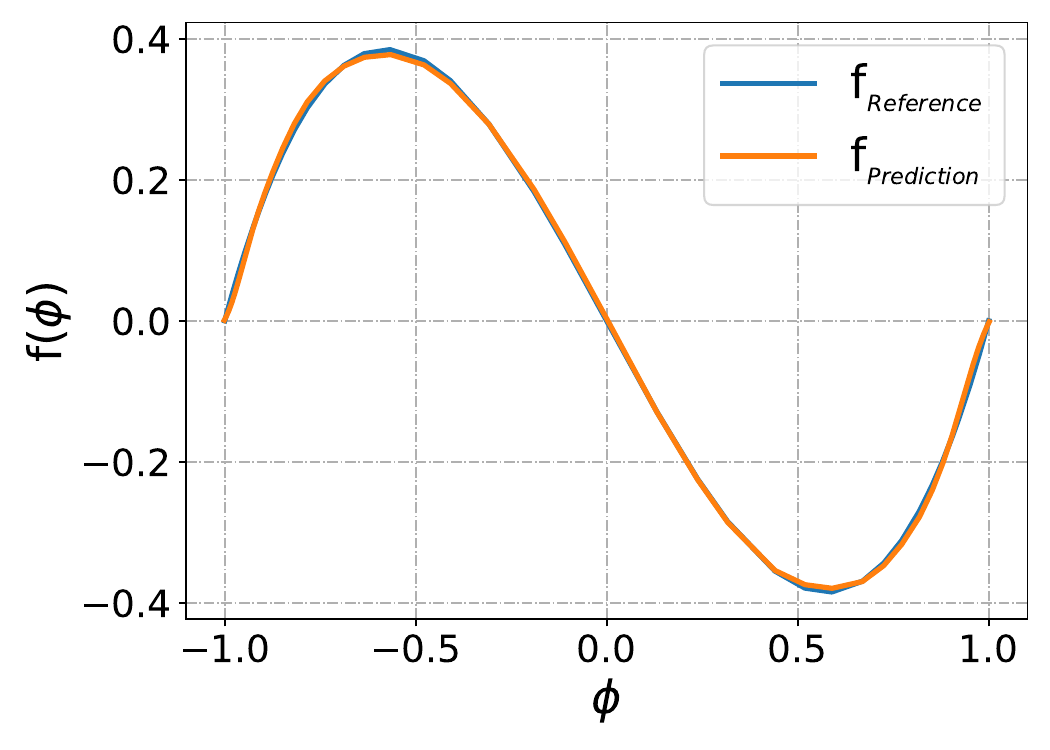}}
    \subfigure[]{\includegraphics[width=0.65\linewidth]{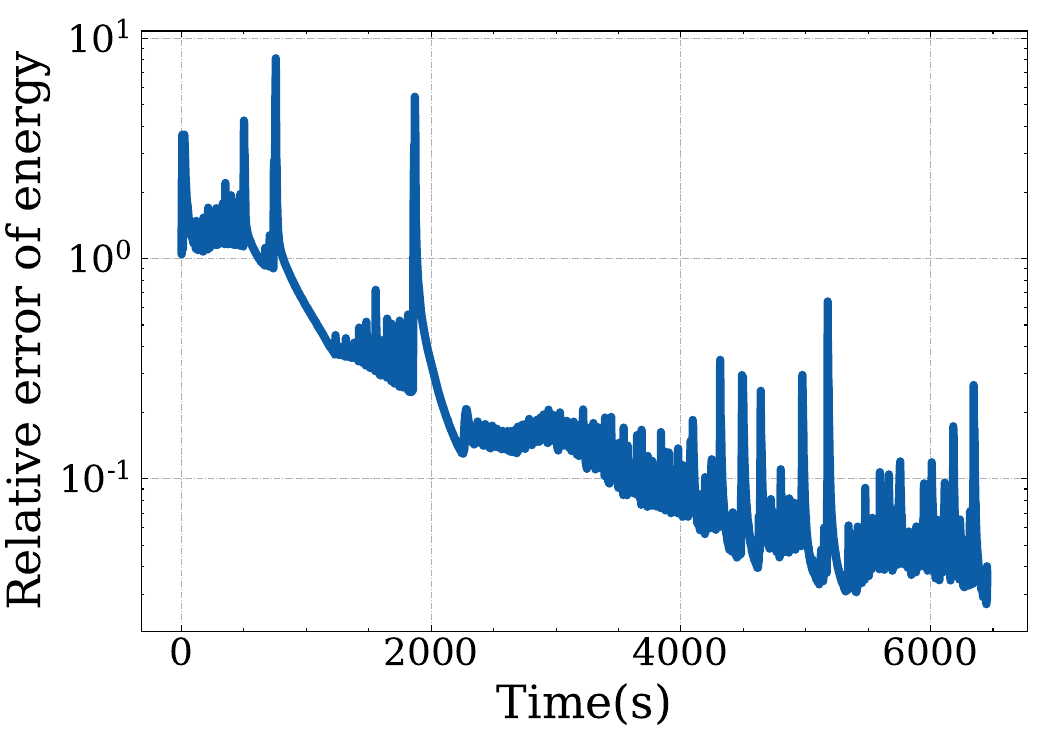}}
    \centering
    \end{minipage}
    \caption{ Results for the inverse problem of 2D Allen-Cahn system. 
    (a)-(b) The reference solutions and the predicted energy functions
    at different time snapshots.
    (c) The plots is the predicted energy function using the PFWNN vs. reference energy function.
    (d) Relative errors for $f^{NN}$ vs. computation times.
     \label{ac2dip_0} }
\end{figure*}

Numerical results of the identifying the energy functional are presented in \figref{ac2dip_0}.
We observe that the predicted energy functional $f$
matches well with the
referenced parameter with high accuracy when the training data is corrupted with noise. 
The predicted energy function $f^{NN}$ is quite close to 
the reference energy function. So the PFWNN can predict the evolution of the energy function accurately. 
Specifically, the resulting relative $L^2$ error in estimating $f(\phi)$ are $1.11e{-3}$.

\begin{figure*}[!ht]
    \centering
    \begin{minipage}{0.6\linewidth}
    \subfigure[]{\includegraphics[width=\linewidth]{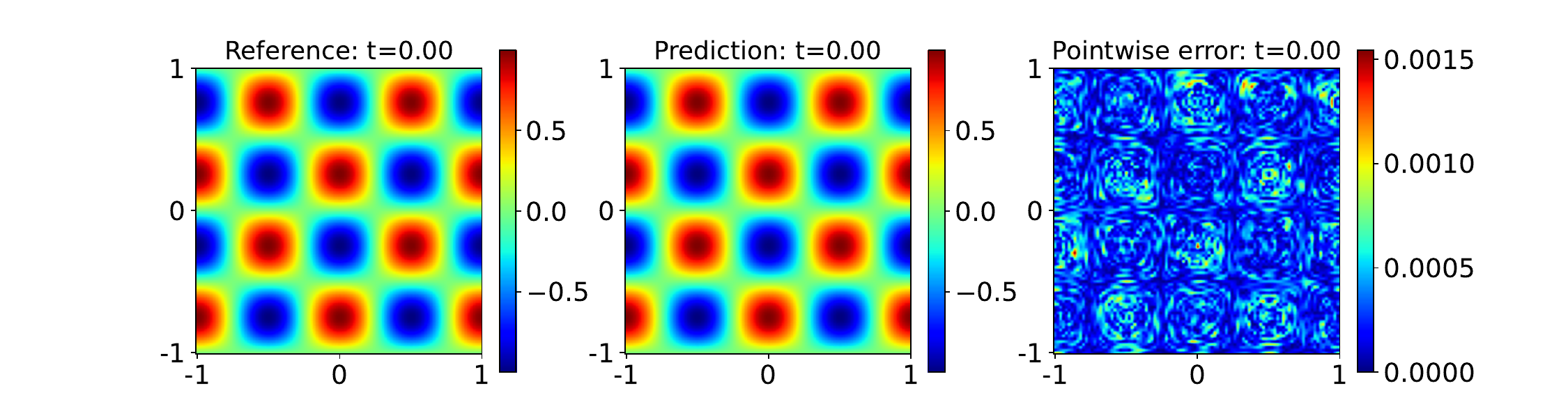}}
    \newline
    \subfigure[]{\includegraphics[width=\linewidth]{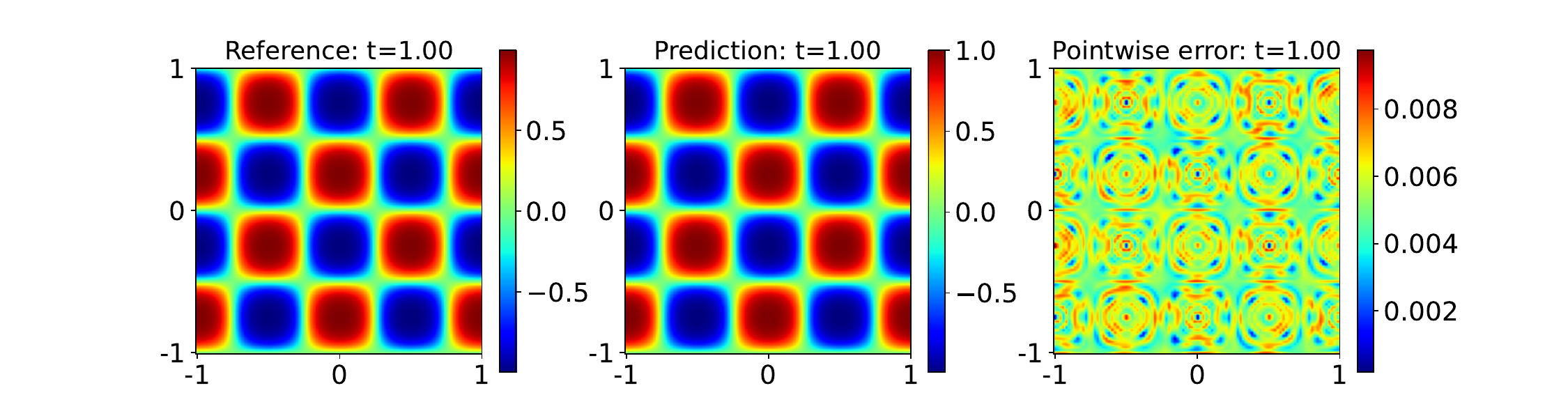}}
    \newline
    \subfigure[]{\includegraphics[width=\linewidth]{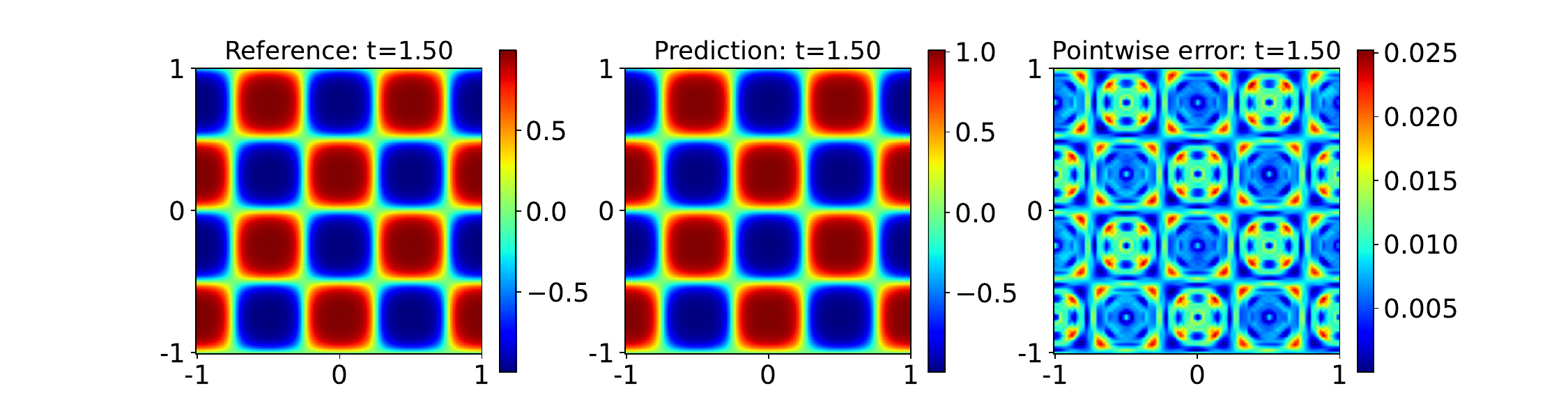}}
    \end{minipage}
    \caption{
    (a)-(c) The reference solutions and the predicted solutions 
    at different time snapshots.
    \label{ac2dip_multi_u} }
\end{figure*}
\begin{figure*}[!ht]
    \centering
    \begin{minipage}{0.60\linewidth}
    \centering
    \subfigure[]{\includegraphics[width=\linewidth]{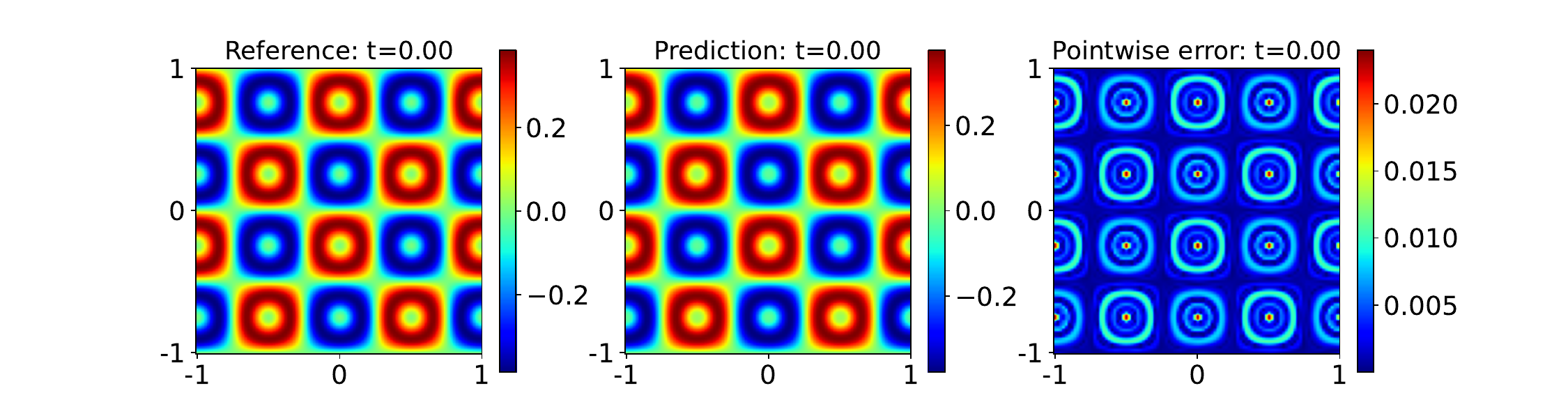}}
    \newline
    \subfigure[]{\includegraphics[width=\linewidth]{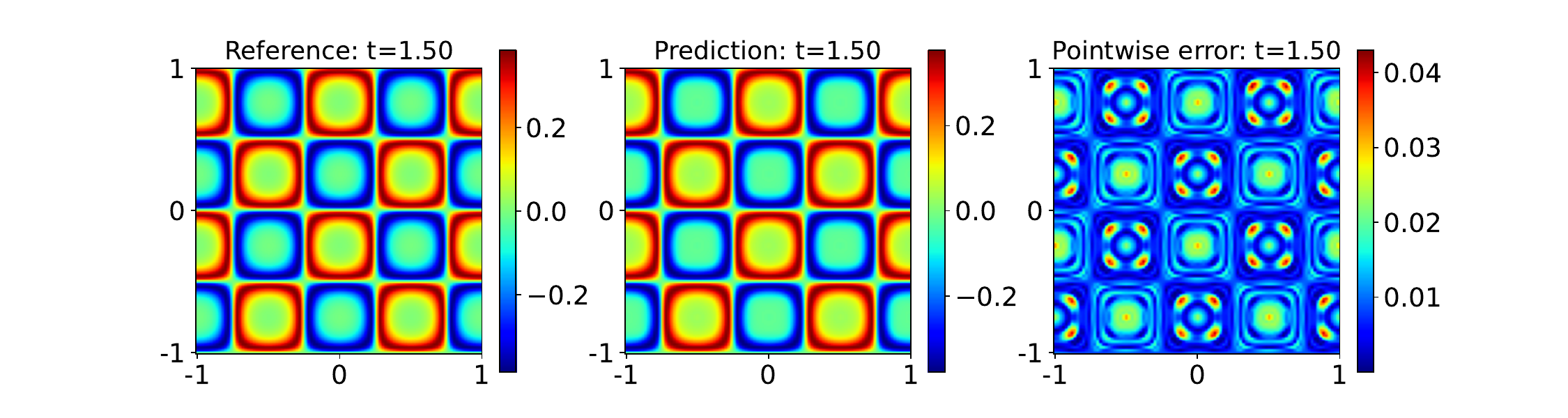}}
    \end{minipage}
    \begin{minipage}{0.35\linewidth}
    \subfigure[]{\includegraphics[width=0.65\linewidth]{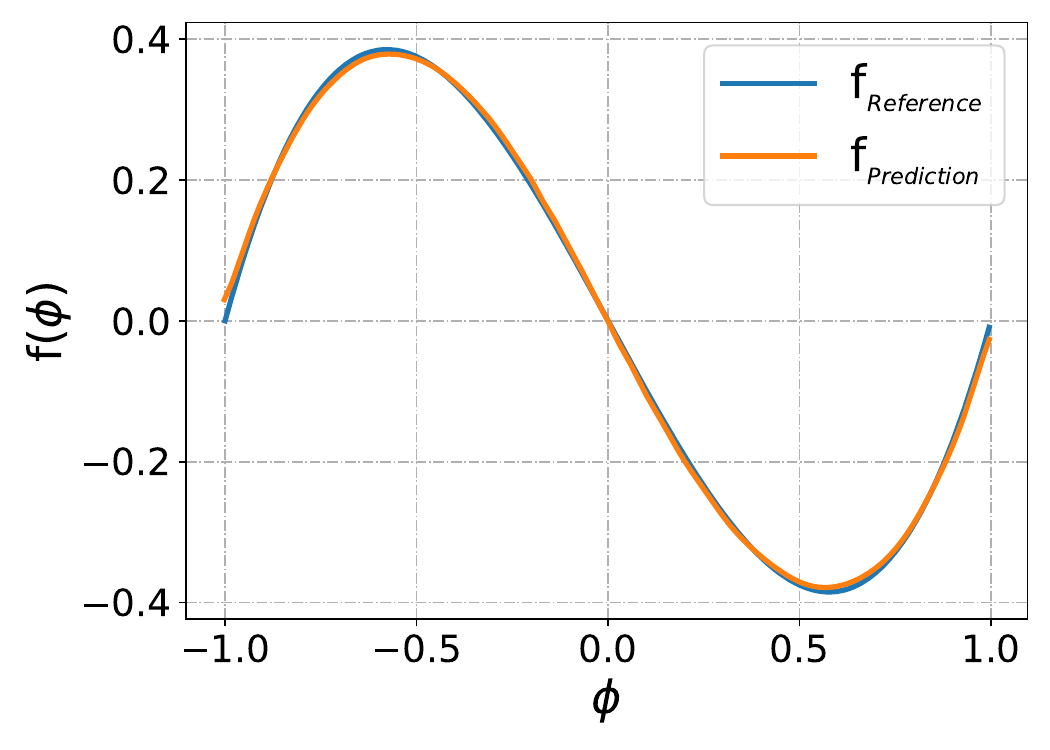}}
    \subfigure[]{\includegraphics[width=0.65\linewidth]{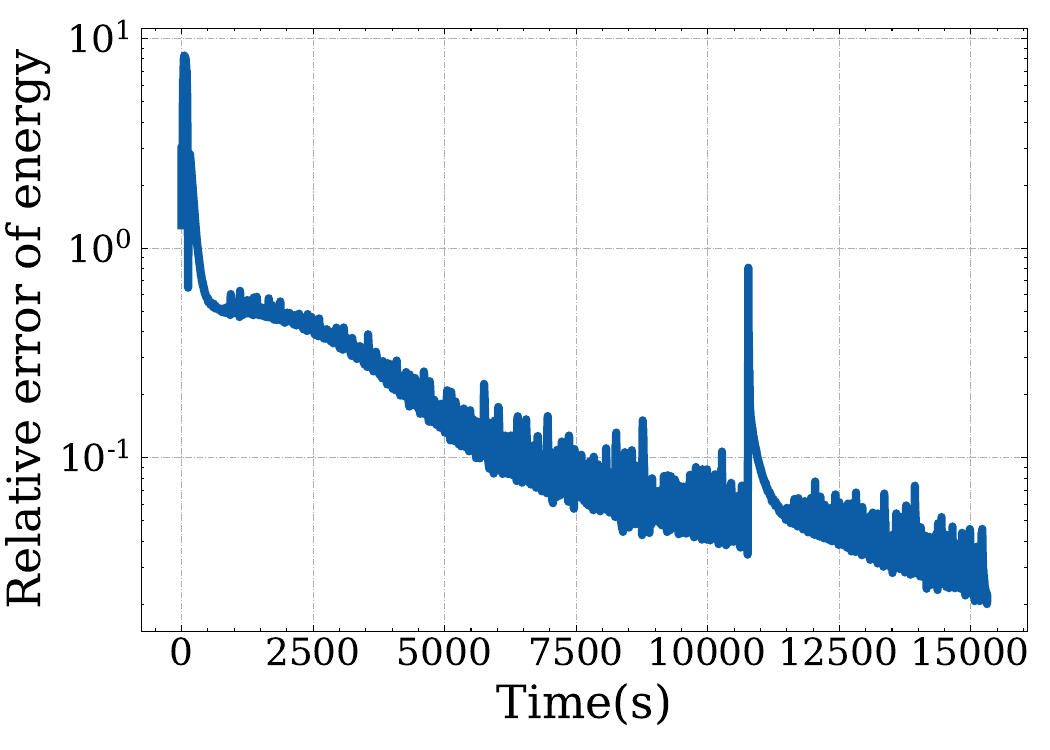}}
    \centering
    \end{minipage}
    \caption{Results for the inverse problem of 2D Allen-Cahn system. 
    (a)-(b) The reference solutions and the predicted energy functions
    at different time snapshots.
    (c) The plots is the predicted energy function using the PFWNN vs. reference energy function.
    (d) Relative errors for $f^{NN}$ vs. computation times.
    \label{ac2dip_multi_energy} }
\end{figure*}

Furthermore, we simulate another example for the Allen-Cahn equation. 
The domain $\Omega=[-1,1]^2$, $ T=2$, and the parameters is set to be
$ M= 1, \lambda^2 = 0.02^2$. 
The initial condition is $\phi(\bm{x}, 0)= sin(2 \pi x_1 )cos(2 \pi x_2)$. 
The neural network architecture and the parameter settings remain the same as 
that of the previous example.
The shape evolution of multiple phase particles is presented in \figref{ac2dip_multi_u}.
The resulting relative $L^2$ error and the maximum absolute error are $1.10e{-3}$ and $3.18e{-3}$, respectively. 
This numerical results also show the ability of the PFWNN to predict
the evolution of the phase-field models. Similarly, numerical results of the inverse problem are presented in \figref{ac2dip_multi_energy}. Specifically, the resulting relative $L^2$ error in estimating $f(\phi)$ are $5.21e{-4}$.

\subsection{Numerical Results for the Cahn-Hilliard Equations}

\subsubsection{The 1D case}
In this section, the Cahn-Hilliard equation 
simulates the phenomenon of phase fusion in 1D situation.
The values of parameters are considered as
$M = 0.01, \lambda^2 = 10^{-4}$. Specifically, we focus on the following specific form,
\begin{equation}\label{eq_ch_0}
    \begin{aligned}
    & \phi_t+\left(  10^{-6} \phi_{x x}- 10^{-2}\left(\phi^3- \phi \right)\right)_{x x}=0, \quad x \in[-1,1], \quad t \in (0,1], \\
    & \phi(0, x)= -\cos (2 \pi x),
    \end{aligned}
\end{equation}

 \begin{figure*}[htb]
    \begin{minipage}{0.8\linewidth}
    \subfigure[]{\includegraphics[width=0.32\linewidth]{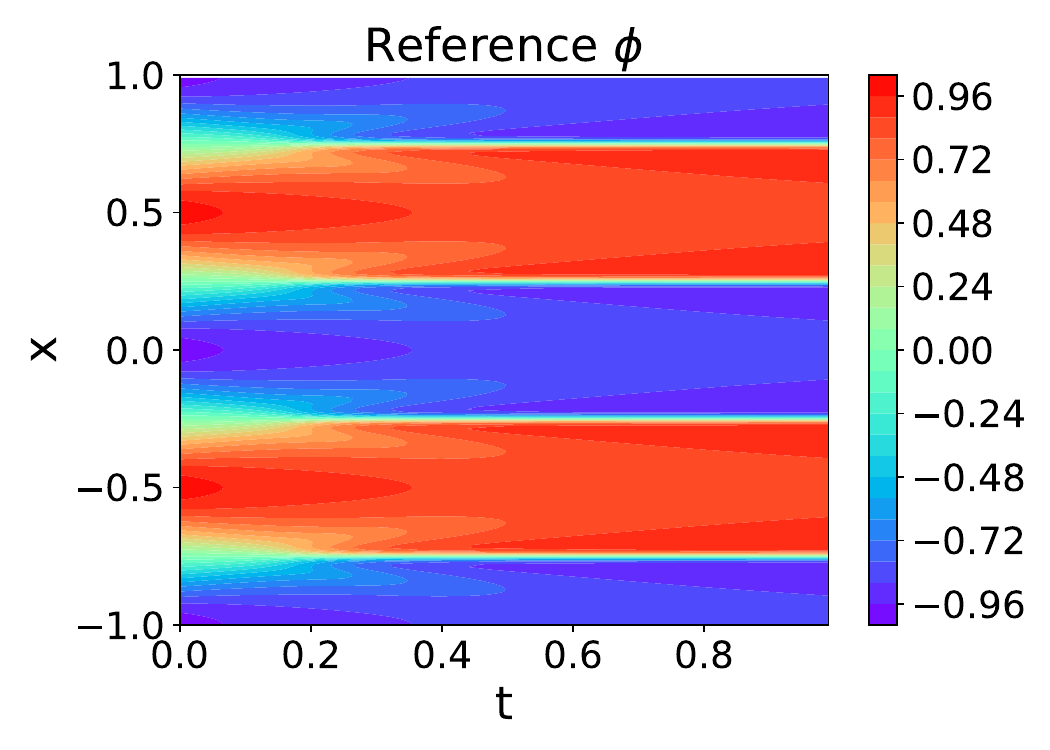}}
    \subfigure[]{\includegraphics[width=0.32\linewidth]{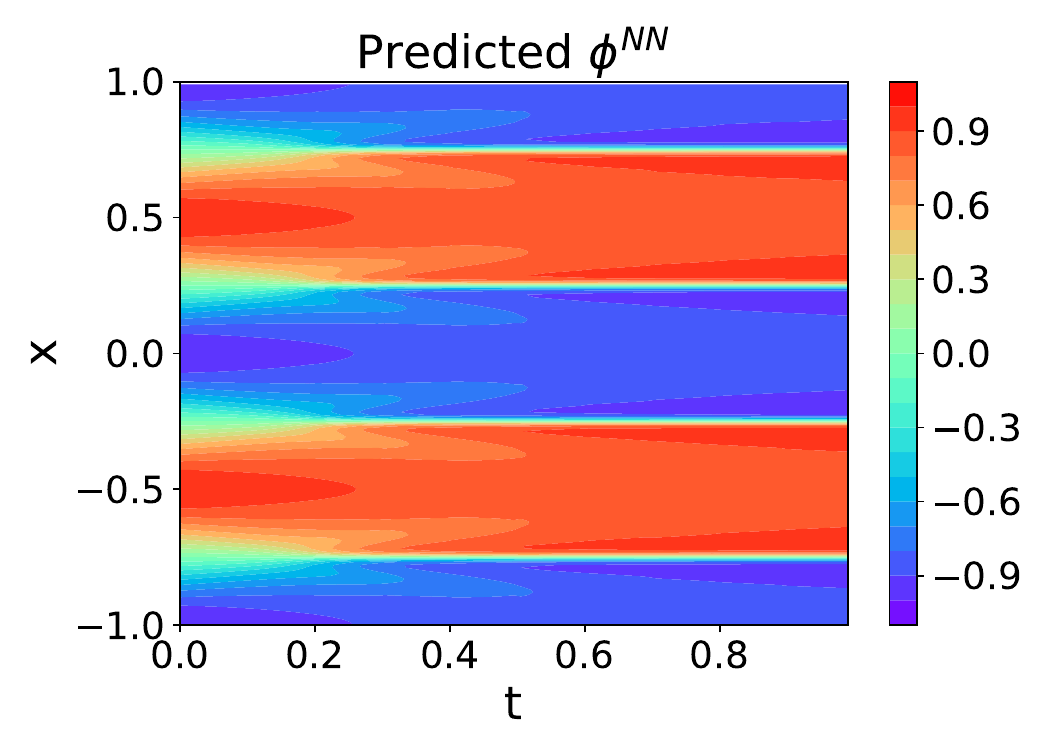}}
    \subfigure[]{\includegraphics[width=0.32\linewidth]{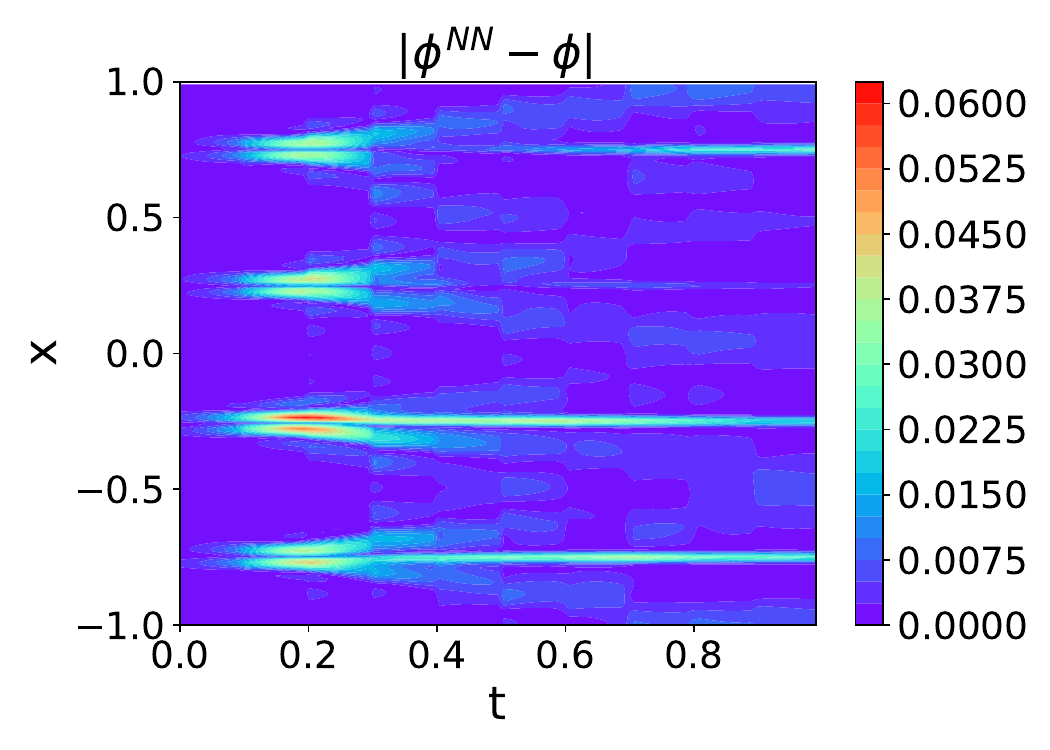}}
    \centering
    \end{minipage}
    \begin{minipage}{0.8\linewidth}
    \subfigure[]{\includegraphics[width=0.32\linewidth]{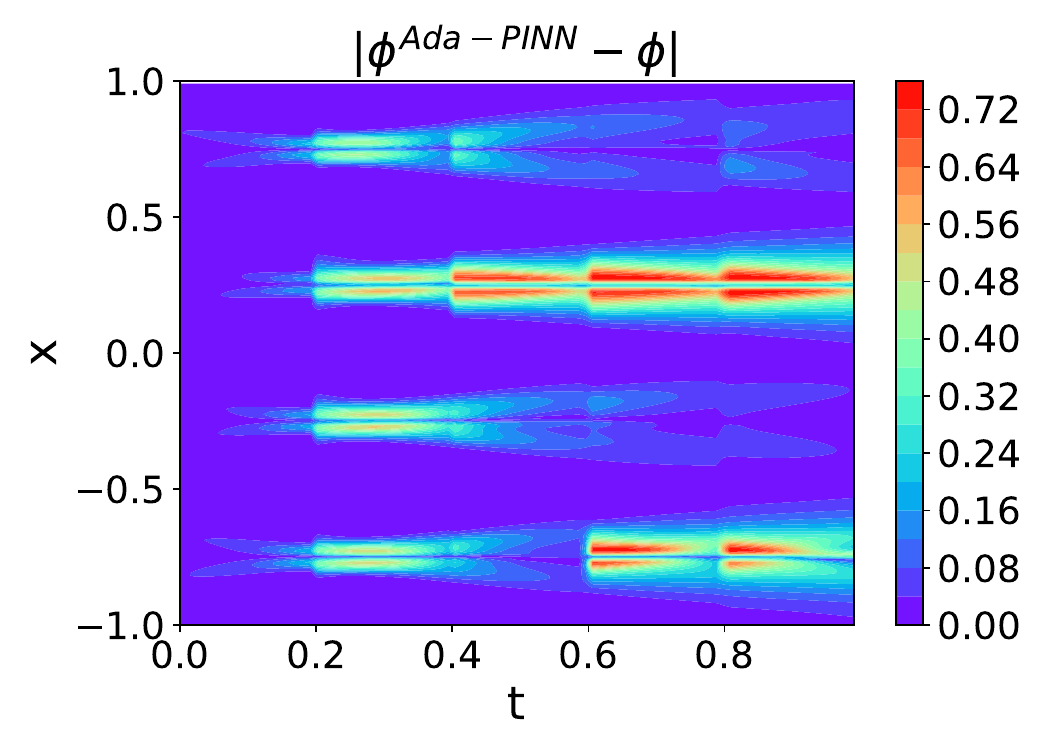}}
    \subfigure[]{\includegraphics[width=0.32\linewidth]{./Figs/AC1d_compare/AC1d_loss_compare.pdf}}
    \subfigure[]{\includegraphics[width=0.32\linewidth]{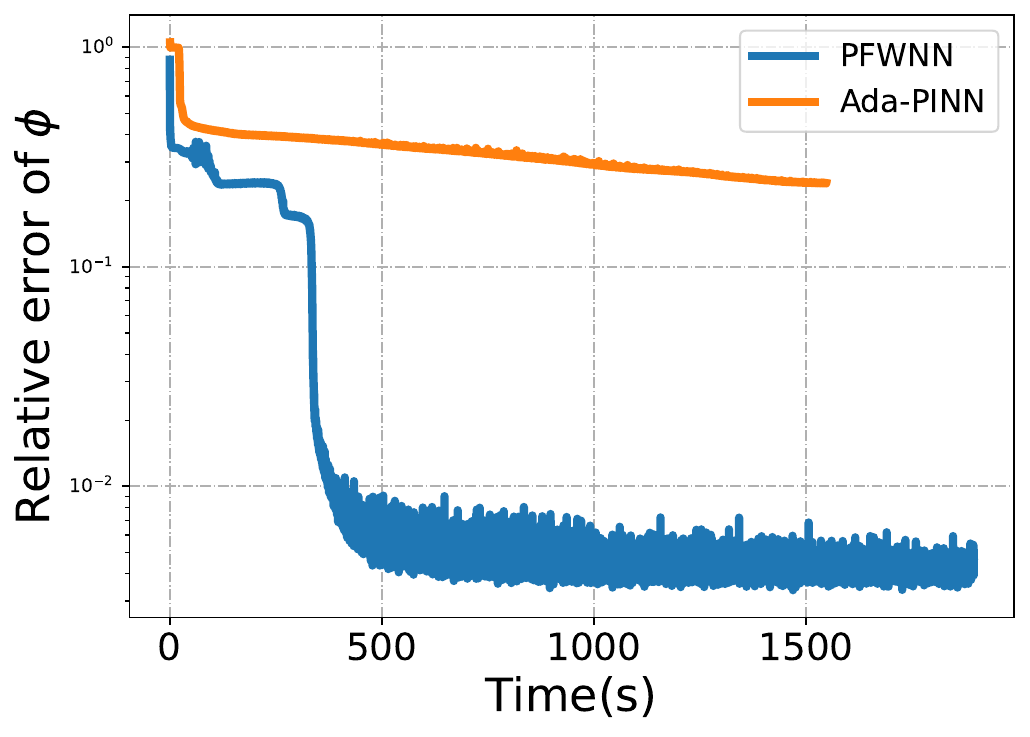}}
    \centering
    \end{minipage}
    \begin{minipage}{0.8\linewidth}
    \subfigure[]{\includegraphics[width=0.7\linewidth]{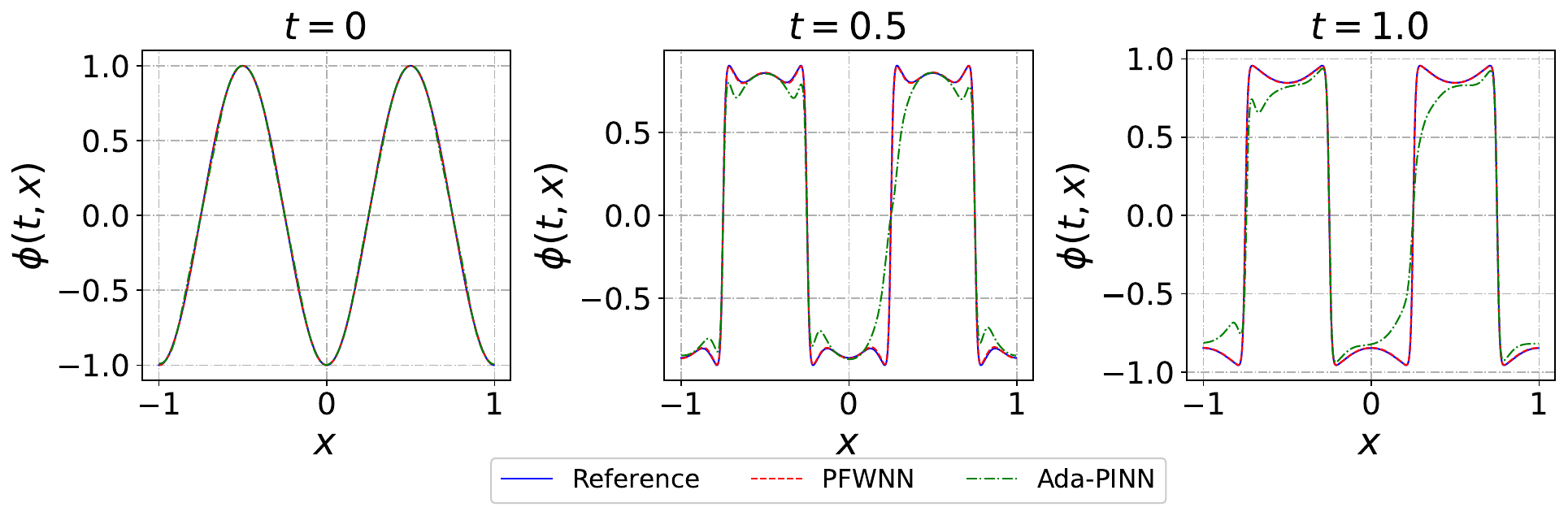}}
    \centering
    \end{minipage}
    \centering
    \caption{Results for 1D Cahn-Hilliard system. Spatio-temporal solutions:
            (a) The reference solution $\phi$, (b) The predicted solution 
            $\phi^{NN}$,
            (c) The point-wise error of the PFWNN.
            (d) The point-wise error of Ada-PINN.
            (e) Loss for solution $\phi^{NN}$ vs. computation times for the last sub-interval.
            (f) Relative errors of solution $\phi^{NN}$ vs. computation times for the last sub-interval.
            (g) The three plots are the predicted solutions of the PFWNN and Ada-PINN vs. reference solutions at different times.\label{ch1dfp}} 
\end{figure*}
Here, the weak form we used is \eqref{CH_mu}. The parameter settings in the experiment are the same as those for the 1D Allen-Cahn equation above, 
but the causality parameter $\epsilon$ is 0.01. 
Beside, compared to the simulation results of the Allen-Cahn equation, 
we divides more timporal intervals to solve the Cahn-Hilliard equation, 
because the phases undergo much greater changes over time.
It is complex and unstable for solving forward problems of the Cahn-Hilliard equation due to the existence of high order derivatives.

The diagrams of the reference solution, the 
predicted solution and the pointwise error 
are shown in \figref{ch1dfp}.
For Ada-PINN, the network does not perfectly learn the sharp curve as the phase field evolves. And the resulting relative $L^2$ error and the maximum absolute error for Ada-PINN are $2.39e{-1}$ and $7.61e{-1}$, respectively.
 The resulting relative $L^2$ error of the PFWNN is $5.44e{-3}$, and the maximum absolute error is $6.26e{-2}$.
It shows that the PFWNN can also accurately predict the solution of the Cahn-Hilliard equations, similar to its performance with the Allen-Cahn equations.

\begin{figure*}[htb]
    \begin{minipage}{0.35\linewidth}
    \subfigure[]{\includegraphics[width=0.9\linewidth]{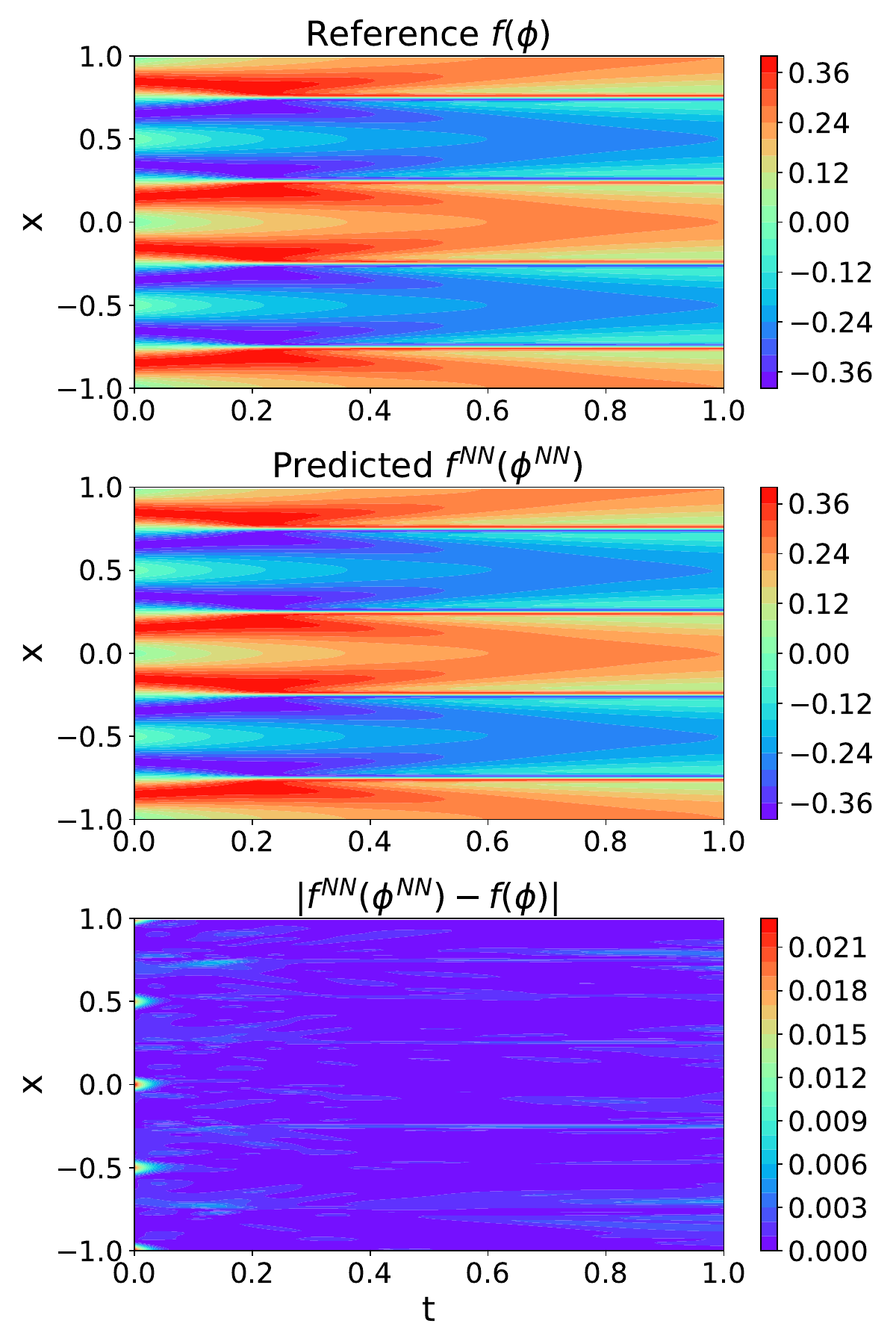}}
    \centering
    \end{minipage}
    \begin{minipage}{0.45\linewidth}
    \subfigure[]{\includegraphics[width=0.65\linewidth]{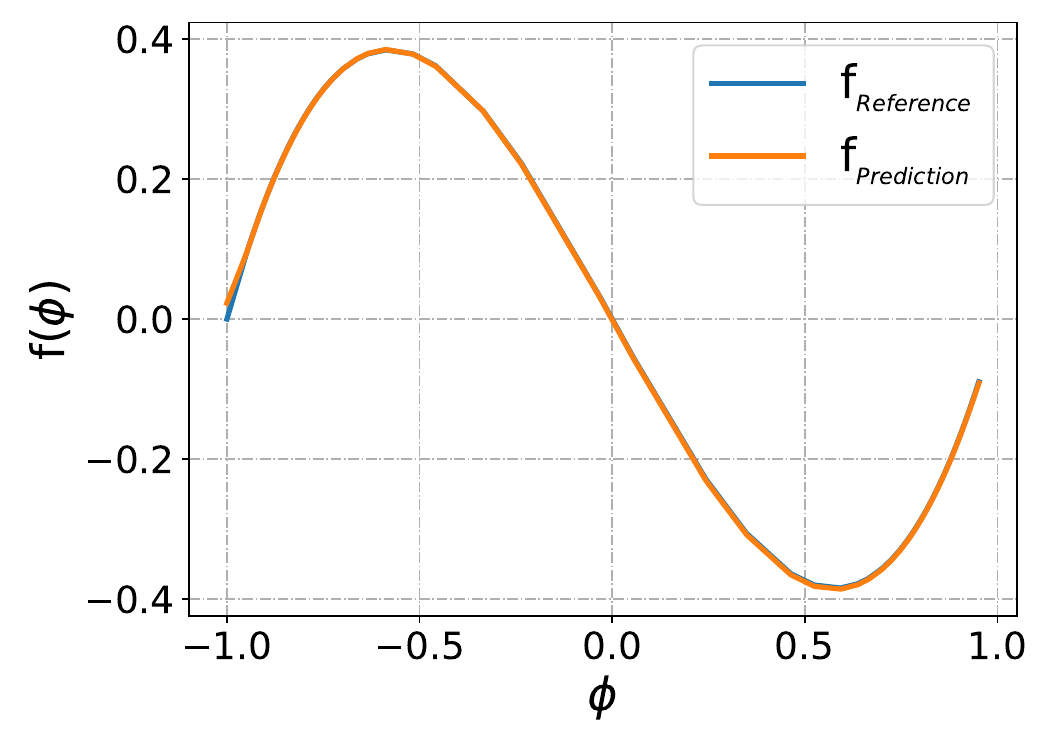}}
    \subfigure[]{\includegraphics[width=0.6\linewidth]{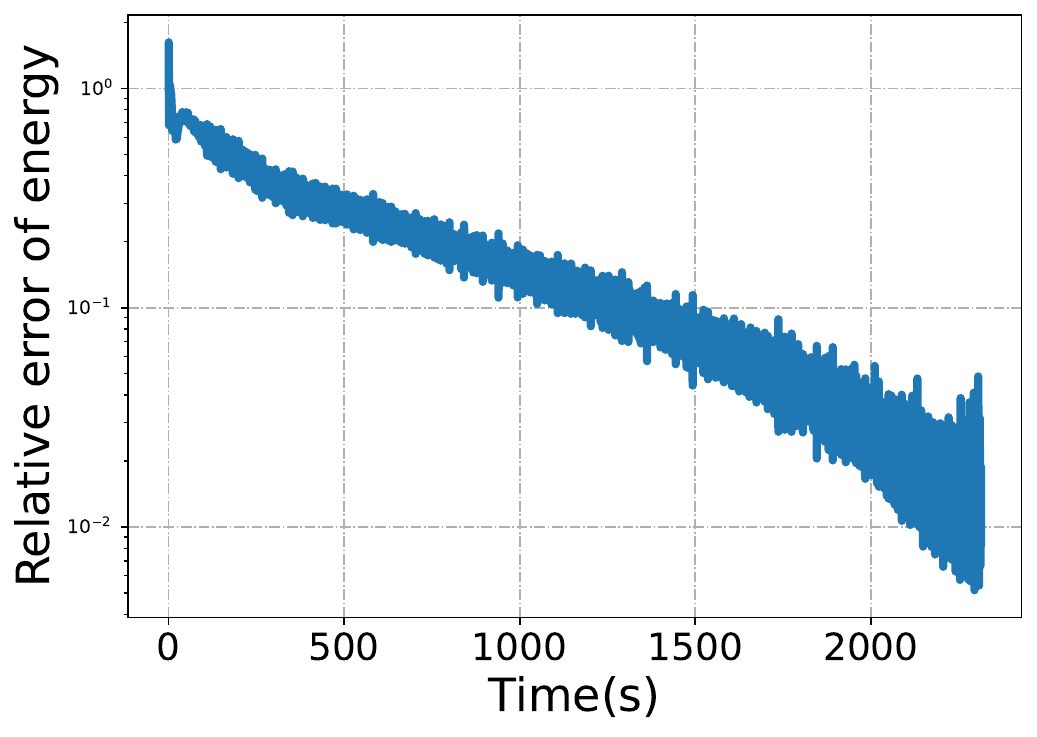}}
    \centering
    \end{minipage}
    \centering
    \caption{Results for the inverse problem of 1D Cahn-Hilliard system. 
            (a) Spatio-temporal energy functional:
            \textit{Top}: The reference energy functional $f(\phi)$, \textit{Middle}: The predicted solution 
            $f^{NN}(\phi^{NN})$,
            \textit{Bottom}: The point-wise error.
            (b) The predicted energy functional $f^{NN}$ vs. reference solutions $f$.
            (c) Relative errors for $f^{NN}$ vs. computation times.
            \label{ch1dip_1d}} 
\end{figure*}

To identify the energy potential of the Cahn-Hilliard equation, we obtain measurement data $\phi(\bm{x}, t_p)$ from the \eqref{eq_ch_0}. Given that the integral region of $f$ is equal to 0, we can shift the identification result to the mean integral.
Numerical results for this example are presented in \figref{ch1dip_1d}. It shows that the PFWNN can correctly identify the unknown parameters $f(\phi)$ with high accuracy.
The predicted energy function $f^{NN}$ is quite close to the reference energy function.
The resulting relative $L^2$ error in estimating $f(\phi)$ of this case is $1.33e{-3}$.
In addition, our advantage in solving the inverse problems is that we can also accurately identify $f$ while the spatio-temporal resolution requirement of the measurement data is lower than that of \cite{brunk2023uniqueness}.

\subsubsection{The 2D case}

Finally, we study the benchmark problem for the Cahn-Hilliard model 
for the two paticles merging. 
The domain is chosen as $\Omega:=[-1,1]^2$, the parameters $M=1, \lambda^2=0.05^2$, and the
initial profile for $\phi$ as
\begin{equation}\label{ch_2d_two}
    \begin{aligned}
    & \phi_t+ \left(0.0025\varDelta^2 \phi-\varDelta \left(\phi^3-\phi\right)\right)=0, \quad x \in[0,1]^2, \quad t \in (0,T], \\
    & \phi(t=0,x,y)=\max \left(tanh \frac{r-R_1}{ \kappa \lambda}, tanh \frac{r-R_2}{ \kappa \lambda}\right),
\end{aligned}
\end{equation}
where $T= 1,r=0.4,\kappa =2$,
$R_1= \sqrt{(x-0.7r)^2+y^2}, R_2= \sqrt{(x+0.7r)^2+y^2}$.
we select the weak form \eqref{CH_mu}, and the weak solution neural network $\phi^{NN}$ with 4 hidden layers 
and 80 nodes in each layer.
During the training process, we 
randomly sample $N_p = 60$ particles in the domain, $N_{int}=15$ integration 
points and $N_{T}=10$ temporal collocation points to evaluate integrals.
Moreover, we set $Iters = 30000$ and $N_{init}=400$ initial collocation points.
The reference data is generated using \eqref{generate}
with $128 \times 128$  modes 
and time-step size of $0.005$. We also obtain measurement data $\phi(\bm{x}, t_p)$ from the reference solution where $\bm{x} \in \Omega ,\ t_p = 0.005p,\  p=0, \ 1, \cdots, \ 100$ for identifying $f$ of \eqref{ch_2d_two}. Then the resulting dataset corresponding to the predicted solution is used for model training, while the remaining data serves as the validation data.

\begin{figure*}[!ht]
    \centering
    \begin{minipage}{0.65\linewidth}
    \subfigure[]{\includegraphics[width=\linewidth]{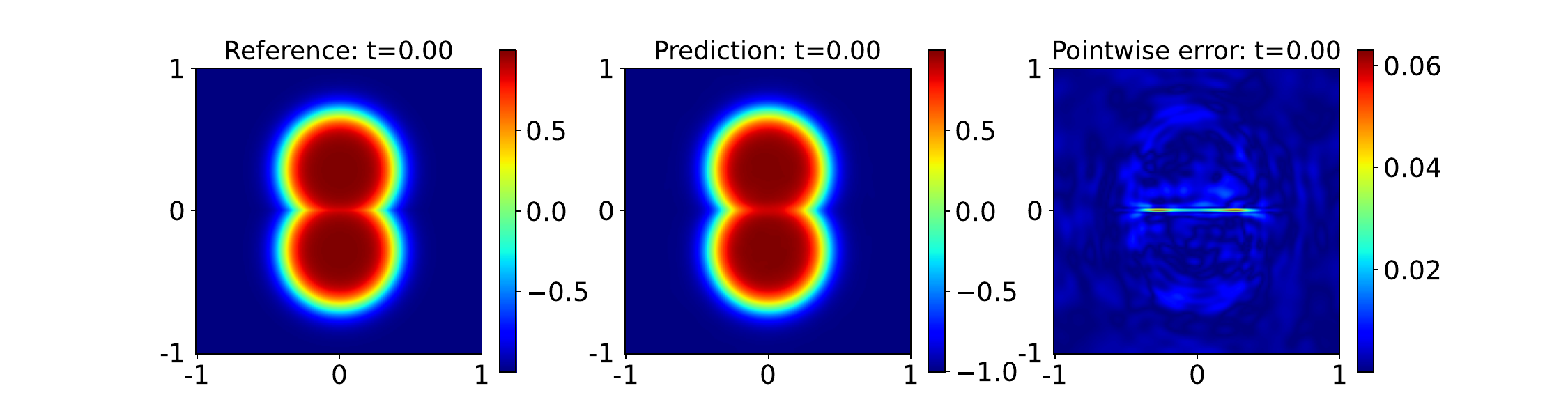}}
    \newline
    \subfigure[]{\includegraphics[width=\linewidth]{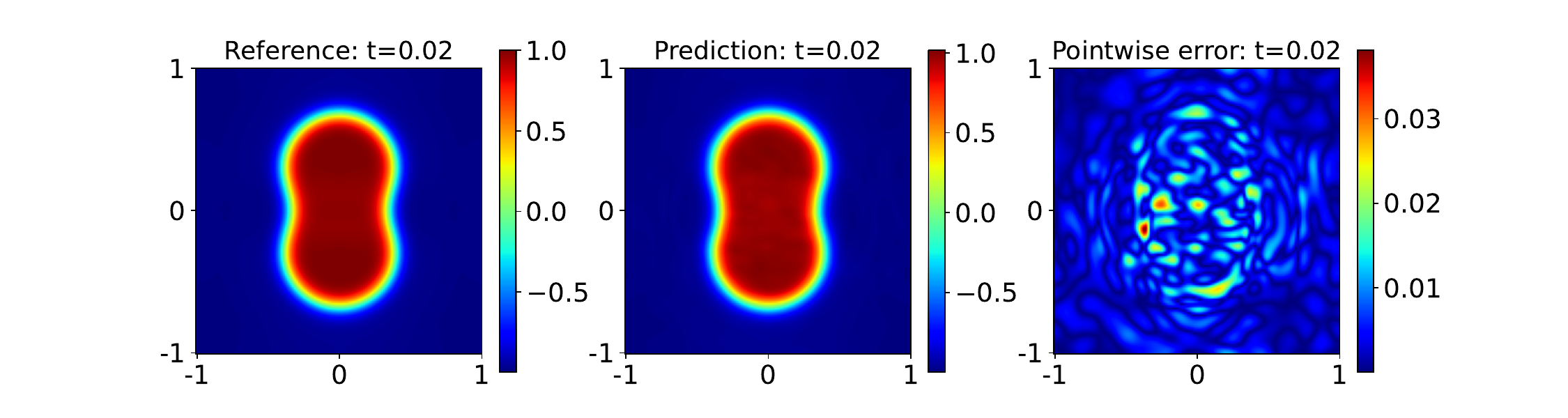}}
    \newline
    \subfigure[]{\includegraphics[width=\linewidth]{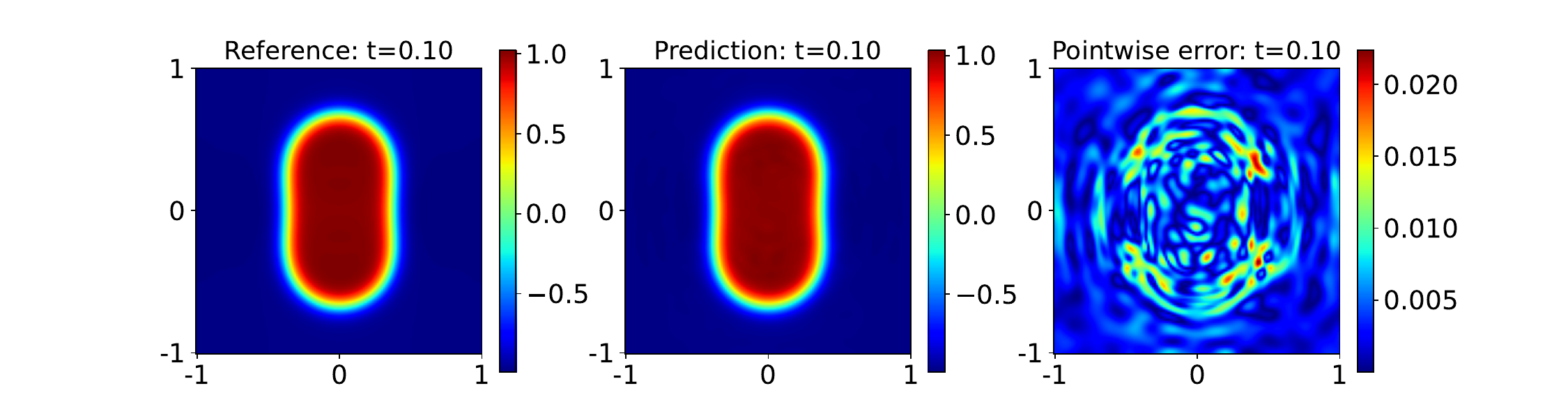}}
    \newline
    \subfigure[]{\includegraphics[width=\linewidth]{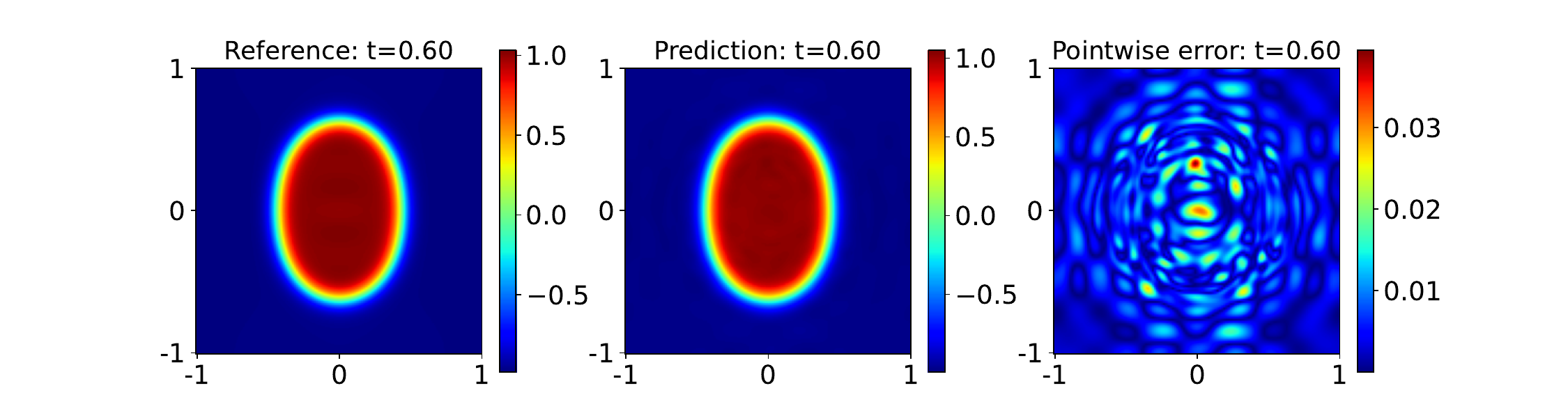}}
    \end{minipage}
    \caption{
    (a)-(d) The reference solutions and the predicted solutions 
    at different time snapshots.
     \label{ch2dip_0} }
\end{figure*}

The temporal evolution of the phase is summarized  
in \figref{ch2dip_0}(a)-(c) for a representative phase-field snapshot 
at time points $t=0.00, 0.02, 0.10, 0.60$.
At the beginning, the two spherical phases diffuse 
from the matrix phase. Subsequently, these two phases 
grow from two spherical to an elliptic. 
The precipitate phases are coarsened and form a 
long strip of tissue along time.
 The resulting relative $L^2$ error and and the maximum absolute error of the PFWNN is $6.57e{-3}$ and $2.68e{-1}$, respectively. 
As a result,
it also provides an accurate approximation of the phase-field diagrams.


\begin{figure*}[!ht]
    \centering
    
    \begin{minipage}{0.6\linewidth}
    \subfigure[]{\includegraphics[width=\linewidth]{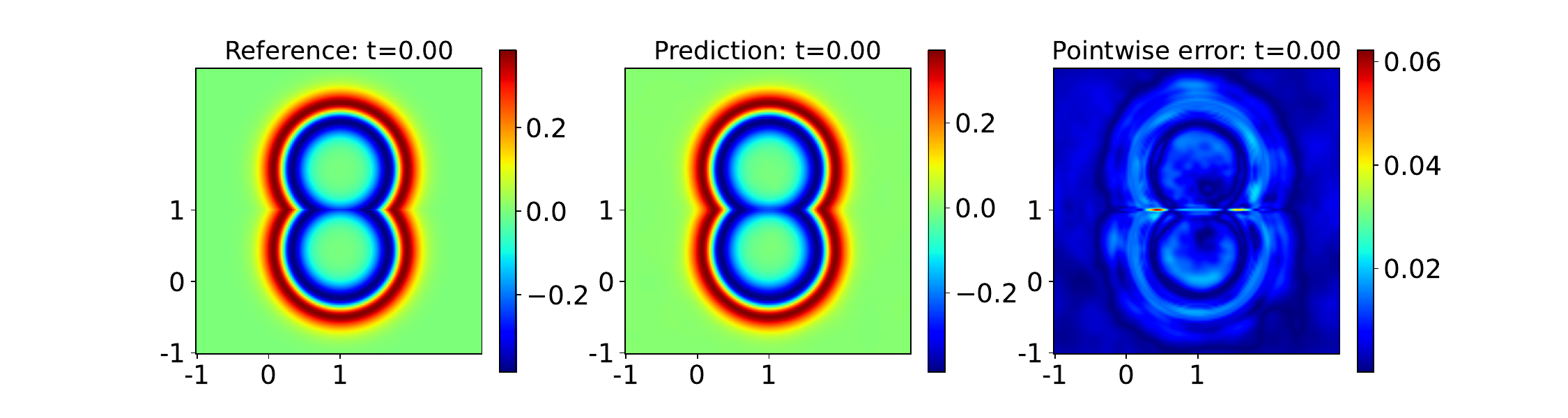}}
    \newline
    \subfigure[]{\includegraphics[width=\linewidth]{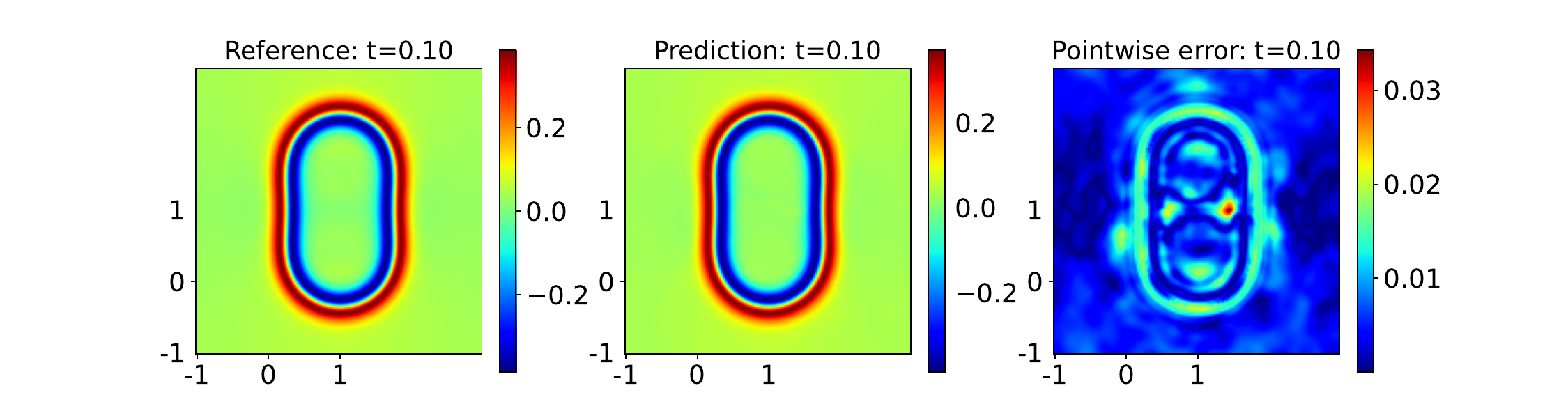}}
    \newline
    \subfigure[]{\includegraphics[width=\linewidth]{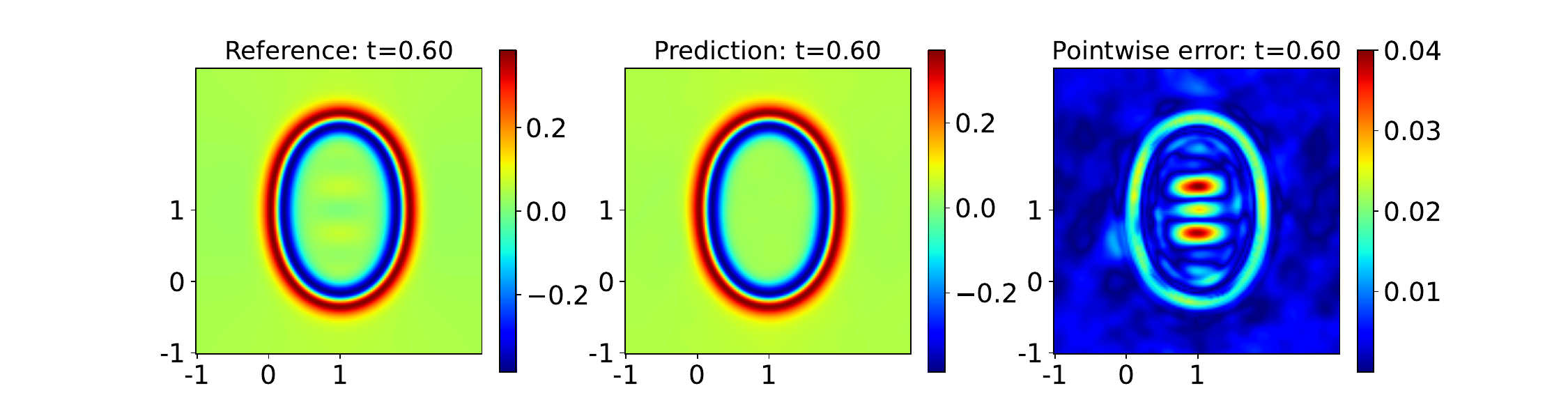}}
    \end{minipage}
    \begin{minipage}{0.35\linewidth}
    \subfigure[]{\includegraphics[width=0.8\linewidth]{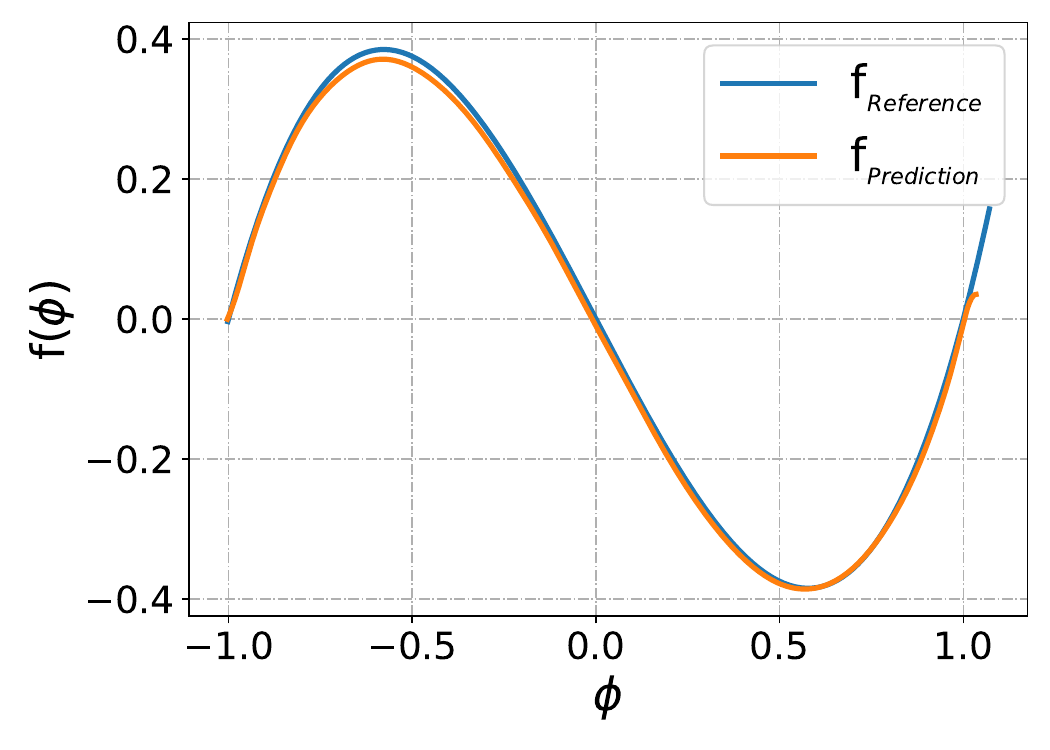}}
    \subfigure[]{\includegraphics[width=0.8\linewidth]{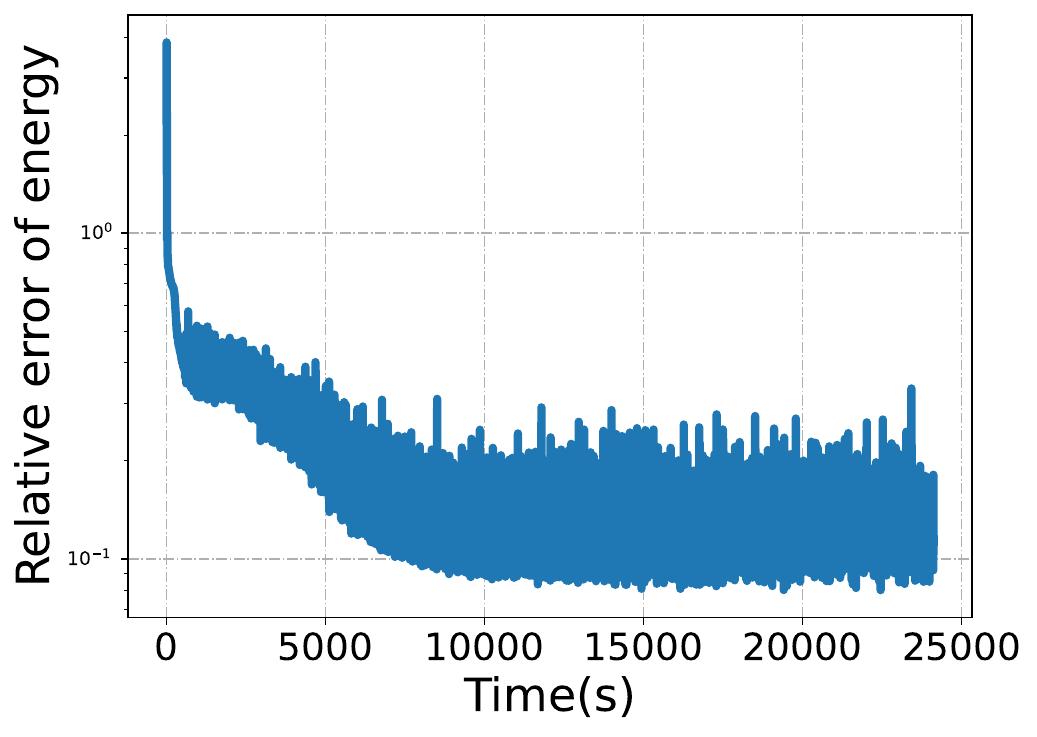}}
    \centering
    \end{minipage}
    
    \caption{Results for the inverse problem of 2D Cahn-Hilliard system. 
    (a)-(c) The reference energy functional and the predicted energy functional
    at different time snapshots.
       (d) The plots is the predicted energy function using the PFWNN vs. reference energy function.
    (e) Relative errors for $f^{NN}$ vs. computation times.
     \label{ch2dfp_0} }
\end{figure*}

The recognition results for the energy functional $f(\phi)$
are shown in \figref{ch2dip_0}(a)-(c).
Although the accuracy of identifying $f(\phi)$ is not as good as that of 
Allen-Cahn equation, 
The predicted energy functional is almost close to the true one.
The resulting relative $L^2$ error in estimating $f(\phi)$ is $6.42e{-3}$.
It can be seen that the 
deep learning method based on weak forms has promising potential in solving the inverse 
problems of high-order and high-dimensional time-dependent PDEs 
from the simulation results of the Cahn-Hilliard equation.

\subsection{The Measurement Data for the Inverse Problems}

The choice of the number of measurement points for the inverse problems plays a crucial role in the accuracy of estimating the energy functional, which involves finding the balance between maximizing accuracy and minimizing computational effort. Consequently, we analyze the effects of the quantity of spatial and temporal samples. Specifically, we employ $N_x$ and $ N_t$ uniform grids in space and time as the measurement data, respectively.
To conduct this investigation, we apply the PFWNN to solve the inverse problems of the 1D Allen-Cahn equation \eqref{eq_ac1d_0}.
\begin{figure*}[!ht]
    \centering
    \includegraphics[width=0.5\textwidth]{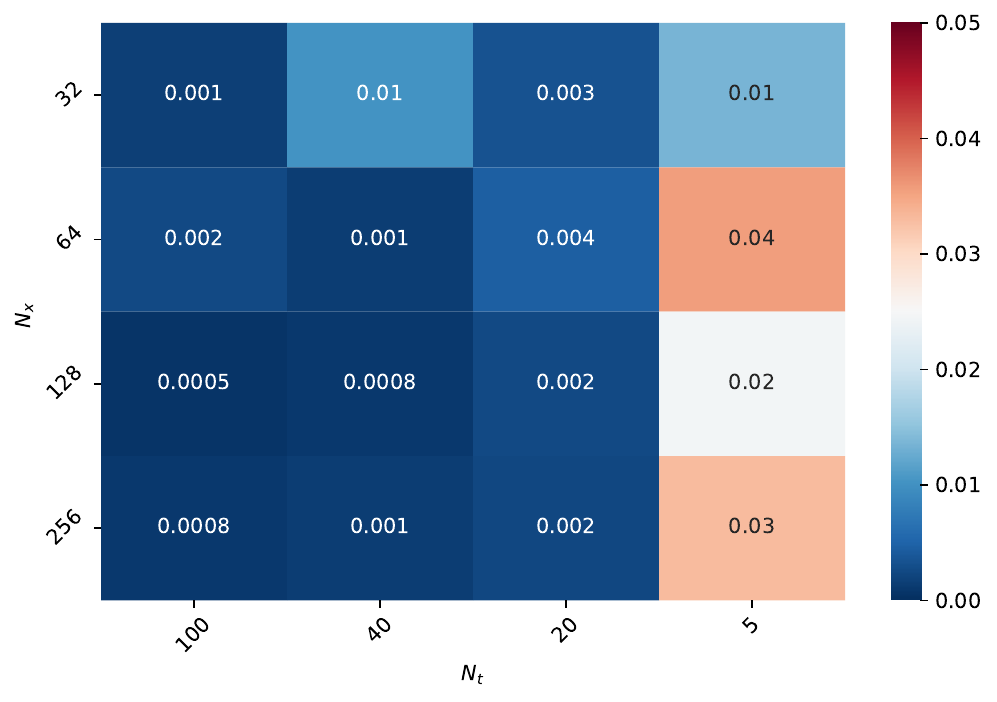}
    \caption{ Performance of the inverse problems for the PFWNN with different combinations of $N_x$ and $N_t$. 
    \label{ac_ch_ip_number} }
\end{figure*}
 Experimental results are presented in \figref{ac_ch_ip_number}. 
 Our ablation study indicates that reducing the number of spatial and temporal data points will decrease accuracy of identifying the energy functional $f$, especially for temporal data points.
In summary, both the number of spatial sampling points and temporal sampling points should be sufficient for effective performance.

\section{Conclusions}
\label{conclusions}

Based on the weak forms, a novel deep learning framework named PFWNN is proposed to 
solve the forward and inverse problems of the Allen-Cahn and Cahn-Hilliard equations.
The original equations are reformulated as their weak forms by integration-by-parts, 
and the high-order equations are transformed into a system of lower-order equations.
The corresponding weak solutions are parameterized as the neural networks, and the test functions are chosen as the locally compactly supported radial basis 
functions in order to render the integration area smaller and more 
concentrated.
As a result, the PFWNN inherits the advantage 
of weak-form methods requiring less regularity of the solution and 
a small number of quadrature points. 
In addition, the periodic boundary conditions are exactly enforced on the weak solutions 
neural networks. The causal training strategy is also adopted to facilitate training of the PFWNN.
For the inverse problems, another neural network is designed for representing the unknown parameters. 
Moreover, both the weak solution neural network
and the network representing the unknown parameters can be effectively trained simultaneously.
We show the efficiency and accuracy of PFWNN and comparison with PINN in several examples and develop the method in detail for 1D and 2D problems of the Allen-Cahn and Cahn-Hilliard equations.
The numerical results confirm that
the PFWNN can work more efficiently with phase-field equations of steep solution and sharp changes,
 and provides an efficient numerical method for the high-dimensional inverse problems.


Despite the excellent performance of the PFWNN for solving 
the Allen-Cahn and Cahn-Hilliard equations, 
there exist some issues such as long-term temporal prediction and non-convexity of the loss function.
Therefore, the future directions involve the development of more efficient time marching 
strategies and loss functions with improved regularity.
Another prospective direction is to solve 
the phase-field models more closely related to physical contexts.

\section{Acknowledgements}
The work was supported in part by National Key Research and Development Program 
of China (No.2023YFA1009100), 
National Natural Science Foundation of China (No.U21A20425), 
and a Key Laboratory of Zhejiang Province.

\bibliographystyle{abbrv}
\bibliography{ref}

\end{document}